\documentclass[a4paper]{article}

\usepackage{amsmath,amssymb,amscd, amsthm}
\newtheorem*{thm}{Theorem}
\newtheorem*{prop}{Proposition}
\newtheorem*{lem}{Lemma}
\newtheorem*{dfn}{Definition}
\newtheorem*{cor}{Corollary}
\newtheorem*{exa}{Example}
\newtheorem*{rem}{Remark}

\newcommand{\Z}{\mathbb{Z}}
\newcommand{\C}{\mathbb{C}}
\newcommand{\R}{\mathbb{R}}
\newcommand{\CP}{\mathbb{P}}

\newcommand{\M}{\mathfrak{M}}
\newcommand{\B}{\mathcal{B}}
\newcommand{\lag}{\mathcal{L}}
\newcommand{\heisen}{\mathcal{H}}
\newcommand{\cli}{\mathcal{C}}
\newcommand{\F}{\mathcal{F}}

\newcommand{\I}{I}

\newcommand{\dimv}{\mathbf{v}}
\newcommand{\dimw}{\mathbf{w}}

\newcommand{\blam}{{\vec{\lambda}}}

\newcommand{\maya}{\mathbf{m}}
\newcommand{\Maya}{\mathcal{M}}

\newcommand{\Hom}{\mathrm{Hom}}
\newcommand{\GL}{\mathrm{GL}}
\newcommand{\gl}{\mathfrak{gl}}
\newcommand{\asl}{\hat{\mathfrak{sl}}}
\newcommand{\ch}{c}
\newcommand{\qf}{\mathcal{R}}
\newcommand{\ec}{H^{\mathrm{mid}}_{S_1}}
\newcommand{\lec}{H^{\mathrm{mid}}_{S_1,\qf}}
\newcommand{\p}{\mathfrak{p}}
\newcommand{\ha}{\mathfrak{h}}

\newcommand{\h}[1]{\frac{#1}{2}}

\newcommand{\quiv}[2]{\M_{\zeta_{#1}}({#2})}

\newcommand{\hilb}{\left(\C^2\right)^{[n]}}

\newcommand{\ueu}{\underline{\mathfrak{e}}}
\newcommand{\eu}{\mathfrak{e}}

\newcommand{\V}{\mathcal{V}}
\newcommand{\tV}{\tilde{\mathcal{V}}}
\newcommand{\tM}{\tilde{M}}

\title{Quiver varieties and Frenkel-Kac construction}
\author{Kentaro Nagao}
\begin{document}

\maketitle

\begin{abstract}
An affine Lie algebra acts on cohomology groups of quiver varieties of affine type.
A Heisenberg algebra acts on cohomology groups of Hilbert schemes of points on a minimal resolution 
of a Kleinian singularity.
We show that in the case of type $A$ the former is obtained by Frenkel-Kac construction from the latter.
\end{abstract}

\tableofcontents

\section{Introduction}

{\ }

\vspace{-12pt}

In this paper we study representation theory associated with the quiver varieties of type $\hat{A}$.

\medskip

A quiver variety is a moduli space of representations of a quiver (\cite{quiver1}). 
Since it is defined as a geometric invariant theory quotient (\cite{GIT}), 
we need to fix a parameter $\zeta$ of the stability condition 
in order to define a quiver variety $\mathcal{M}_\zeta$. 
The space of parameters of stability conditions has a chamber structure. 
A wall consists of parameters where the stability condition and the semistability condition are not equivalent. 
When the parameter cross a wall, the variety is changed by a flop, as is typical in the geometric invariant theory (\cite{GIT-flip}).
But in our situation it is also known that the underlying $C^\infty$-manifold is not changed (\cite{quiver1}). 
In particular, there exists a canonical isomorphism between cohomology groups of quiver varieties associated with different chambers.

\medskip

One can construct a level-$1$ integrable highest weight representation of a Kac-Moody Lie algebra on a direct sum of cohomology groups of quiver varieties 
associated with a specific parameter $\zeta_0$ (\cite{quiver2}).

On the other hand, one can also construct the Fock space representation of a Heisenberg algebra on a direct sum of cohomology groups 
of Hilbert schemes of points on a surface (\cite{nakajima-heisen}).
When the surface is a minimal resolution of a Kleinian singularity, 
the Hilbert schemes of points on it can be described as quiver varieties of affine type, which are associated with another specific parameter $\zeta_\infty$.

An affine Lie algebra contains the Heisenberg algebra as its subalgebra. 
One can reconstruct a level-$1$ integrable representation of an affine Lie algebra 
from a level-$1$ representation of the Heisenberg algebra (\textbf{Frenkel-Kac construction} \cite{frenkel-kac}). 

So both cohomology groups of $\mathcal{M}_{\zeta_0}$ and $\mathcal{M}_{\zeta_\infty}$ are endowed with 
level-$1$ integrable highest weight actions of an affine Lie algebra. 
Our purpose is to show that the canonical isomorphism of cohomology groups induced by the diffeomorphism intertwines those actions.

\medskip

In this paper we will deal with the case of type $A$ only. Let us explain our method. 

The quiver varieties are endowed with specific $S^1$-actions. 
The quiver varieties associated with different chambers are $S^1$-equivariantly diffeomorphic. 
The actions of the affine Lie algebra and the Heisenberg algebra on the cohomology groups can be lifted to the $S^1$-equivariant cohomology groups. 
We will show the isomorphism induced by the $S^1$-equivariant diffeomorphism intertwine the two actions on the $S^1$-equivariant cohomology groups, 
which will be followed by the result for the ordinary cohomology groups. 

It is an advantage that the $S^1$-equivariant cohomology groups have specific bases indexed by the $S^1$-fixed points. 
Moreover the $S^1$-fixed points of $\mathcal{M}_{\zeta_0}$ and $\mathcal{M}_{\zeta_\infty}$ can be parametrized in term of Young diagrams. 
We will exhibit the two actions of the affine Lie algebra with respect to these bases 
and show that the operation called "taking cores and quotients of Young diagrams" induces an isomorphism as representations.  

It seems to have been known or believed that the fixed points of the quiver varieties associated with different stability conditions are related  by the operation taking cores and quotients. 
For example, in \S 7.2.4 of \cite{haiman-lec-note} it is tried to characterize the wreath Macdonald polynomials in term of cores and quotients, 
where the wreath Macdonald polynomials are conjectured to exist, to correspond to fixed points of quiver varieties, and so, to depend on the chamber of the stability conditions. 
We will prove that this relation can be realized by the diffeomorphism. 
We should mention that this result is already used in \S 6.3 of \cite{framed-moduli-affine-Lie-alg}, 
but neither proof nor reference is given there.     

\medskip

The paper is organised as follows.
In \S \ref{pre} we are devoted to combinatorial preliminaries. 
In \S \ref{Frenkel-Kac construction} we review Frenkel-Kac construction and give an explicit formula for the resulting representation.
In \S \ref{Quiver varieties} we study the geometry of the quiver varieties and show the correspondence of the fixed points is given by the operation taking cores and quotients.
In \S \ref{Geometric representation assosiated with quiver varieties} we review Nakajima's geometrical construction of the representations of the affine Lie algebra and the Heisenberg algebra, 
and describe the representations on the $S^1$-equivariant cohomology groups explicitly using the bases derived from the fixed points.
The results of \S \ref{Frenkel-Kac construction}--\S \ref{heisenberg} are combined to give the theorem for the equivariant cohomology groups in \S \ref{main for ec}. 
In \S \ref{ordinary} we deduce the result for the ordinary cohomology groups.

\subsection*{Acknowledgements}
This is a part of the master thesis written under the supervision of Professor Hiraku Nakajima. 
The author would like to thank him for his valuable comments, warm encouragement and careful proofreading.
The author also would like to thank him writing the paper \cite{nakajima-moduli-on-ALE}.

\section{Preliminaries}\label{pre}

\subsection{Notations for indices}

Fix an integer $l$ larger than $2$.
In this paper we sometimes use half integers in order to describe everything symmetrically.
Symbols $k$ and $h$ are used for half integers.

We set $I=\{0,\ldots,l-1\}$.  
We sometimes identify $I$ with the set of modulo $l$ equivalent classes of integers. 
We also identify $I$ with the set of vertices of the cyclic quiver.

We set $\tilde{I}=\{\h{1},\h{3},\ldots,l-\h{1}\}$. 
We sometimes identify $\tilde{I}$ with the set of modulo $l$ equivalent classes of half integers. 
We also identify $\tilde{I}$ with the set of edges of the cyclic quiver.

\subsection{Notations for Young diagrams}

Let $\Pi$ denote the set of all Young diagrams.
Identify a Young diagram with a subset of $(\Z_{\geq0})^2$.
A {\bf node} is an element of $(\Z_{\geq 0})^2$.
The {\bf content} of a node $(a,b)$ is the number $a-b$.
A node is called {\bf $i$-node} if its content equals to $i$ modulo $l$. 
For $\lambda\in\Pi$ we define
\[
n_j(\lambda)=\sharp\{(a,b)\in\lambda\mid a-b=j\},\quad
f_\lambda(z)=\sum_{j\in\Z}n_j(\lambda)z^j\in\Z[z^\pm]
\]
and
\[
v_{{i}}(\lambda)=\sharp\{(a,b)\in\lambda\mid a-b\equiv i\ (\mathrm{mod}\,l)\},\quad
\dimv(\lambda)=(v_{{i}}(\lambda))\in\Z^I.
\]

For $\lambda\in\Pi$, a node $(a,b)$ is called {\bf addable} 
if $(a,b)\notin\lambda$ and $(a-1,b),(a,b-1)\in\lambda$.
A node $(a,b)$ is called {\bf removable} 
if $(a,b)\in\lambda$ and $(a+1,b),(a,b+1)\notin\lambda$.
Let $A_{\lambda,i}$ (resp. $R_{\lambda,i}$) denote the set of all addable (removable) $i$-nodes for $\lambda$,

Ler $X$ be a removable or an addable node of $\lambda$. We set 
\begin{align*}
\eta^+(\lambda,i,X)=&\,\sharp\{\text{addable $i$-node to the right of $X$}\}\\
&-\sharp\{\text{removable $i$-node to the right of $X$}\},\\
\eta^-(\lambda,i,X)=&\,\sharp\{\text{addable $i$-node to the left of $X$}\}\\
&-\sharp\{\text{removable $i$-node to the left of $X$}\}.
\end{align*}

\subsection{Notations for Maya diagrams}

A {\bf Maya diagram} is 
an increasing sequence of half integers $\maya=(k_j)_{j\geq 1}$ such that $k_{j+1}=k_j+1$ for sufficiently large $j$.
Let $\Maya$ denote the set of all Maya diagrams.

Note that a Maya diagram can be identified with a map $\Z+\h{1}\to\{\pm 1\}$ such that 
\[
\maya(h)=
\begin{cases}
1 & \quad \text{for} \ h\gg 0\\
-1& \quad \text{for} \ h\ll 0.
\end{cases}
\]
We define the {\bf charge} of $\maya$ by
\[
c(\maya)=\sharp\{h<0\mid \maya(h)=1\}-\sharp\{h>0\mid \maya(h)=-1\}.
\]
For $c\in\Z$, we define $\maya^{(c)}$ by
$\maya^{(c)}(h)=\maya(h+c)$.
Note that $c(\maya^{(c)})=c(\maya)-c$.

\subsection{Young diagrams and Maya diagrams}

\subsubsection{}\label{prebosonfermionmap}

For $\lambda\in\Pi$, let us define $\maya_{\lambda}\in\Maya$.
First, note that
\[
n_{k-\h{1}}(\lambda)-n_{k+\h{1}}(\lambda)=
\begin{cases}
-1\ \,\text{or}\ \,0 & (k<0),\\
0\ \,\text{or}\ \,1 & (k>0).
\end{cases} 
\]
Then we define
\[
\maya_{\lambda}(k)=
\begin{cases}
1 & (k<0,\ n_{k-\h{1}}(\lambda)-n_{k+\h{1}}(\lambda)=-1),\\
-1 & (k<0,\ n_{k-\h{1}}(\lambda)-n_{k+\h{1}}(\lambda)=0),\\
1 & (k>0,\ n_{k-\h{1}}(\lambda)-n_{k+\h{1}}(\lambda)=0),\\
-1 & (k>0,\ n_{k-\h{1}}(\lambda)-n_{k+\h{1}}(\lambda)=1).
\end{cases} 
\]
\begin{center}
%WinTpicVersion3.08
\unitlength 0.1in
\begin{picture}( 28.0000, 16.6500)( 12.0000,-24.6500)
% LINE 2 0 3 0
% 4 2600 2000 1400 800 2600 2000 3800 800
% 
\special{pn 8}%
\special{pa 2600 2000}%
\special{pa 1400 800}%
\special{fp}%
\special{pa 2600 2000}%
\special{pa 3800 800}%
\special{fp}%
% LINE 0 0 3 0
% 2 1800 1200 2000 1000
% 
\special{pn 20}%
\special{pa 1800 1200}%
\special{pa 2000 1000}%
\special{fp}%
% LINE 2 0 3 0
% 2 2000 1000 2200 1200
% 
\special{pn 8}%
\special{pa 2000 1000}%
\special{pa 2200 1200}%
\special{fp}%
% LINE 2 0 3 0
% 2 2200 1200 2400 1400
% 
\special{pn 8}%
\special{pa 2200 1200}%
\special{pa 2400 1400}%
\special{fp}%
% LINE 0 0 3 0
% 2 2400 1400 2600 1200
% 
\special{pn 20}%
\special{pa 2400 1400}%
\special{pa 2600 1200}%
\special{fp}%
% LINE 0 0 3 0
% 2 2600 1200 2800 1000
% 
\special{pn 20}%
\special{pa 2600 1200}%
\special{pa 2800 1000}%
\special{fp}%
% LINE 0 0 3 0
% 2 2800 1000 3000 800
% 
\special{pn 20}%
\special{pa 2800 1000}%
\special{pa 3000 800}%
\special{fp}%
% LINE 2 0 3 0
% 2 3000 800 3200 1000
% 
\special{pn 8}%
\special{pa 3000 800}%
\special{pa 3200 1000}%
\special{fp}%
% LINE 2 0 3 0
% 2 3200 1000 3400 1200
% 
\special{pn 8}%
\special{pa 3200 1000}%
\special{pa 3400 1200}%
\special{fp}%
% LINE 0 0 3 0
% 2 3400 1200 3800 800
% 
\special{pn 20}%
\special{pa 3400 1200}%
\special{pa 3800 800}%
\special{fp}%
% LINE 2 0 3 0
% 2 1200 2400 4000 2400
% 
\special{pn 8}%
\special{pa 1200 2400}%
\special{pa 4000 2400}%
\special{fp}%
% DOT 2 0 3 0
% 13 1500 2400 1700 2400 1900 2400 2100 2400 2300 2400 2500 2400 2700 2400 2900 2400 3100 2400 3300 2400 3500 2400 3700 2400 3700 2400
% 
\special{pn 8}%
\special{sh 1}%
\special{ar 1500 2400 10 10 0  6.28318530717959E+0000}%
\special{sh 1}%
\special{ar 1700 2400 10 10 0  6.28318530717959E+0000}%
\special{sh 1}%
\special{ar 1900 2400 10 10 0  6.28318530717959E+0000}%
\special{sh 1}%
\special{ar 2100 2400 10 10 0  6.28318530717959E+0000}%
\special{sh 1}%
\special{ar 2300 2400 10 10 0  6.28318530717959E+0000}%
\special{sh 1}%
\special{ar 2500 2400 10 10 0  6.28318530717959E+0000}%
\special{sh 1}%
\special{ar 2700 2400 10 10 0  6.28318530717959E+0000}%
\special{sh 1}%
\special{ar 2900 2400 10 10 0  6.28318530717959E+0000}%
\special{sh 1}%
\special{ar 3100 2400 10 10 0  6.28318530717959E+0000}%
\special{sh 1}%
\special{ar 3300 2400 10 10 0  6.28318530717959E+0000}%
\special{sh 1}%
\special{ar 3500 2400 10 10 0  6.28318530717959E+0000}%
\special{sh 1}%
\special{ar 3700 2400 10 10 0  6.28318530717959E+0000}%
\special{sh 1}%
\special{ar 3700 2400 10 10 0  6.28318530717959E+0000}%
% LINE 2 2 3 0
% 2 1900 2400 1900 1100
% 
\special{pn 8}%
\special{pa 1900 2400}%
\special{pa 1900 1100}%
\special{dt 0.045}%
% LINE 2 2 3 0
% 4 2500 2400 2500 1300 2500 1300 2500 1300
% 
\special{pn 8}%
\special{pa 2500 2400}%
\special{pa 2500 1300}%
\special{dt 0.045}%
\special{pa 2500 1300}%
\special{pa 2500 1300}%
\special{dt 0.045}%
% LINE 2 2 3 0
% 4 2700 2400 2700 1300 2700 1300 2700 1100
% 
\special{pn 8}%
\special{pa 2700 2400}%
\special{pa 2700 1300}%
\special{dt 0.045}%
\special{pa 2700 1300}%
\special{pa 2700 1100}%
\special{dt 0.045}%
% LINE 2 2 3 0
% 4 2900 2400 2900 1300 2900 1300 2900 900
% 
\special{pn 8}%
\special{pa 2900 2400}%
\special{pa 2900 1300}%
\special{dt 0.045}%
\special{pa 2900 1300}%
\special{pa 2900 900}%
\special{dt 0.045}%
% LINE 2 2 3 0
% 4 3700 2400 3700 1300 3700 1300 3700 900
% 
\special{pn 8}%
\special{pa 3700 2400}%
\special{pa 3700 1300}%
\special{dt 0.045}%
\special{pa 3700 1300}%
\special{pa 3700 900}%
\special{dt 0.045}%
% LINE 2 2 3 0
% 4 3500 2400 3500 1300 3500 1300 3500 1100
% 
\special{pn 8}%
\special{pa 3500 2400}%
\special{pa 3500 1300}%
\special{dt 0.045}%
\special{pa 3500 1300}%
\special{pa 3500 1100}%
\special{dt 0.045}%
% STR 2 0 3 0
% 3 1900 2450 1900 2550 5 0
% $-\h{7}$
\put(19.0000,-25.5000){\makebox(0,0){$-\h{7}$}}%
% STR 2 0 3 0
% 3 2500 2450 2500 2550 5 0
% $-\h{1}$
\put(25.0000,-25.5000){\makebox(0,0){$-\h{1}$}}%
% STR 2 0 3 0
% 3 2700 2450 2700 2550 5 0
% $\h{1}$
\put(27.0000,-25.5000){\makebox(0,0){$\h{1}$}}%
% STR 2 0 3 0
% 3 2900 2450 2900 2550 5 0
% $\h{3}$
\put(29.0000,-25.5000){\makebox(0,0){$\h{3}$}}%
% STR 2 0 3 0
% 3 3500 2450 3500 2550 5 0
% $\h{9}$
\put(35.0000,-25.5000){\makebox(0,0){$\h{9}$}}%
% STR 2 0 3 0
% 3 3700 2450 3700 2550 5 0
% $\h{11}$
\put(37.0000,-25.5000){\makebox(0,0){$\h{11}$}}%
% STR 2 0 3 0
% 3 3900 2450 3900 2550 5 0
% $\cdots$
\put(39.0000,-25.5000){\makebox(0,0){$\cdots$}}%
\end{picture}%
\end{center}
\vspace{1mm}

\subsubsection{}\label{eta}

An integer $j$ such that $\left(\maya_\lambda(j-\h{1}), \maya_\lambda(j+\h{1})\right)=(1,-1)$ (resp. $(-1,1)$)
corresponds to an addable (resp. a removable) node of $\lambda$. 
Its content equals to $j$. 

For an addable or a removable node with content $j$ we have
\begin{align*}
&\eta^-(X,j,\lambda)\\
&=\,\sharp\left\{j'\equiv j\ \Big|\ j'<j,\ \left(\maya_\lambda\left(j-\h{1}\right), \maya_\lambda\left(j'+\h{1}\right)\right)=(1,-1)\right\}\\
&\quad -\sharp\left\{j'\equiv j\ \Big|\ j'<j,\ \left(\maya_\lambda\left(j-\h{1}\right), \maya_\lambda\left(j'+\h{1}\right)\right)=(-1,1)\right\}\\
&=\,\sharp\left\{h\equiv j-\h{1}\ \Big|\ h<j,\ \maya_\lambda(h)=1\right\}-\sharp\left\{h\equiv j+\h{1}\ \Big|\ h<j,\ \maya_\lambda(h)=1\right\}.
\end{align*}

\subsubsection{}\label{bosonfermionmap}

For $\lambda\in\Pi$ we have $c(\maya_{\lambda})=0$.
Conversely given $\maya\in\Maya$ with charge $0$, 
there exists a unique Young diagram $\lambda$ such that $\maya_{\lambda}=\maya$.
Consequently, we have the bijection
\[
\begin{array}{cccc}
F\colon&\Z\times\Pi&\longrightarrow&\Maya\\ 
&(c,\lambda)&\longmapsto&\ \ \maya_{\lambda}^{(-c)}.
\end{array}
\]
Let us define $q(\maya)\in\Pi$ by $F^{-1}(\maya)=(\ch(\maya),q(\maya))$.

\subsection{Cores and quotients}

\subsubsection{}\label{corequotient}

For $k\in \tilde{I}$ and $\maya\in\Maya$, 
we define $\maya_k\in\Maya$ by
\[
\maya_k(h)=\maya\left(l\left(h-\h{1}\,\right)+k\right).
\]
Note that $\maya$ can be recovered from $\{\maya_k\}_{k\in\tilde{I}}$. 
We have $\ch(\maya)=\sum\ch(\maya_k)$.

For $\lambda\in\Pi$, we set $c_k(\lambda)=c(\maya_{\lambda,k})$ and $q_k(\lambda)=q(\maya_{\lambda,k})$.
We define the {\bf $l$-core} of $\lambda$ by $\mathbf{c}(\lambda)=(c_k(\lambda))\in(\Z^{\tilde{I}})_0=\{(c_\h{1},\ldots,c_{l-\h{1}})\in\Z^{\tilde{I}}\mid \sum{c_k}=0\}$. 
%where $\left(\Z^I\right)_0$ is the sublattice of $\Z^I$ defined by the equation that the sum of the coordinate function equals to $0$.
We also define the {\bf $l$-quotient} of $\lambda$ by $\mathbf{q}(\lambda)=q_k(\lambda)\in\Pi^{\tilde{I}}$.
We have the bijection
\[
CQ=\mathbf{c}\times \mathbf{q}\colon \Pi\longrightarrow(\Z^{\tilde{I}})_0\times\Pi^{\tilde{I}}.
\]

\subsubsection{}\label{corequotientrem}
\begin{lem}
\[
c_k(\lambda)=v_{k-\h{1}}(\lambda)-v_{k+\h{1}}(\lambda).
\]
\end{lem}
\begin{proof}
Note that for For $k\in{\tilde{I}}$ and $h\in\Z+\h{1}$ the sign of $h$ and of $l\left(h-\h{1}\,\right)+k$ coincide.
\begin{align*}
c_k(\lambda)&=c(\maya_{\lambda,k})\\
&=\sharp\{h<0\mid\maya_{\lambda,k}(h)=1\}-\sharp\{h>0\mid\maya_{\lambda,k}(h)=-1\}\\
&=\sharp\left\{l\left(h-\h{1}\,\right)+k<0\ \Big|\ \maya_{\lambda}\left(l\left(h-\h{1}\,\right)+k\right)=1\right\}\\
&\quad-\sharp\left\{l\left(h-\h{1}\,\right)+k>0\ \Big|\ \maya_{\lambda}\left(l\left(h-\h{1}\,\right)+k\right)=-1\right\}\\
&=\sum_{h<0}\left(n_{l\left(h-\h{1}\,\right)+k-\h{1}}(\lambda)-n_{l\left(h-\h{1}\,\right)+k+\h{1}}(\lambda)\right)\\
&\quad-\sum_{h>0}\left(-n_{l\left(h-\h{1}\,\right)+k-\h{1}}(\lambda)+n_{l\left(h-\h{1}\,\right)+k+\h{1}}(\lambda)\right)\\
&=v_{k-\h{1}}(\lambda)-v_{k+\h{1}}(\lambda).
\end{align*}
\end{proof}

\subsubsection{}\label{cqprop}
\begin{prop}
For $\lambda\in\Pi$ we set $\underline{\lambda}=CQ^{-1}(\mathbf{c}(\lambda),(\emptyset,\ldots,\emptyset))$. Then we have
\[
f_\lambda(z)=f_{\underline{\lambda}}(z)+\sum_kz^{lc_k(\lambda)}f_{q_k(\lambda)}(z^l)\left(z^{k-l+\h{1}}+\cdots+z^{k-\h{1}}\right).
\]
\end{prop}
\begin{proof}
Let $\lambda^+$ be a Young diagram obtained from $\lambda$ by adding a node with content $j\in\Z$ to $q_k(\lambda)$.
Since $k_{q_k(\lambda)}(h)=k_\lambda\left(l(c_k(\lambda)+h-\h{1})+k\right)$ we have 
\begin{align*}
&k_\lambda\left(l(c_k(\lambda)+j-1)+k\right)=1,\\
&k_\lambda\left(l(c_k(\lambda)+j)+k\right)=-1,\\
&k_{\lambda^+}\left(l(c_k(\lambda)+j-1)+k\right)=-1,\\
&k_{\lambda^+}\left(l(c_k(\lambda)+j)+k\right)=1,
\end{align*}
and $k_\lambda(h)=k_{\lambda^+}(h)$ if $h\neq l(c_k({\lambda})+j-1)+k,l(c_k({\lambda})+j)+k$.
This means 
\[
n_{j'}(\lambda^+)-n_{j'}(\lambda)=\begin{cases}1&l(c_k(\lambda)+j-1)+k+\h{1}\leq j'\leq l(c_k(\lambda)+j)+k-\h{1},\\0&\text{otherwise}.\end{cases}
\]
Then the statement follows by induction on $\sum|q_k(\lambda)|$.
\end{proof}

\section{Frenkel-Kac construction}\label{Frenkel-Kac construction}
{\ }

\vspace{-14pt}

In this section we review Frenkel-Kac construction.

Let $Q$ denote the finite root lattice of type $A_{l-1}$ and 
$\heisen_{l-1}=\heisen_{\asl_{l}}$ denote the Heisenberg algebra of type $A_{l-1}$.
Let $\B^{\otimes \tilde{I}}$ be a vector space with a basis indexed by $\Pi^{\tilde{I}}$. 
This is endowed with an action of $\heisen_{l-1}$.
By Frenkel-Kac construction, we get an action $\asl_l$ on $\C Q\otimes \B^{\otimes \tilde{I}}$.
By the operation taking $l$-cores and $l$-quotients, 
$Q\times \Pi^{\tilde{I}}$ is bijective to $\Pi$. 
Thus $\C Q\otimes \B^{\otimes \tilde{I}}$ has a basis indexed by $\Pi$.
We will give an explicit formula of the action of $\asl_{l}$ with respect to this basis (Theorem \ref{formulaforFrenkelKac}).
This formula seems to be well-known to many people, but I can not find any references.
We will prove the formula by reducing to the \textbf{boson-fermion correspondence}.

\subsection{Boson-fermion correspondence}

\subsubsection{}\label{bosonicfockspace}

The \textbf{Heisenberg algebra} $\heisen$ is the $\C$-algebra generated by $p(m),p(-m)\ (m\in\Z_{>0})$ and the central element $c$ with relation
\[
[p(m),p(n)]=\delta_{0,n+m}c.
\]
The Heisenberg algebra $\heisen$ acts on the polynomial ring $\C[p_1,p_2,\ldots]$ with infinitely many variables by
\[
p(-m)\cdot f=mp_mf,\quad p(m)\cdot f=\frac{\partial}{\partial p_m}f,\quad c\cdot f=f.
\]

Identify $p_m$ with the power sum symmetric function, then $\C[p_1,p_2,\ldots]$ is isomorphic to the 
ring of symmetric functions with infinitely many variables.  
Identify the Schur function $s_{\lambda}$ with a formal element $b_{\lambda}$, then the 
ring of symmetric function with infinitely many variables is isomorphic to $\B=\oplus_{\lambda\in\Pi}\,\C b_{\lambda}$ as vector space.

Hence we have the action of $\heisen$ on $\B$. This representation is called the {\bf bosonic Fock space}.

\subsubsection{}\label{Fermionicfockspace}

The \textbf{Clliford algebra} $\cli$ is the $\C$-algebra generated by $\psi(k),\psi^*(k)\ (k\in\Z+\h{1})$ and the central element $c$ with relation
\[
\{\psi(k),\psi(h)\}=\{\psi^*(k),\psi^*(h)\}=0,\ \{\psi(k),\psi^*(h)\}=\delta_{k,h}c.
\]

The Clliford algebra $\cli$ acts on $\F=\oplus_{\maya\in\Maya}\,\C\maya$ by
\begin{align*}
\psi_k(\maya)&=
\begin{cases}
(-1)^{i-1}(\ldots,k_{i-1},k_{i+1},\ldots) & \quad \text{if}\  k_i=k \ \text{for some i},\\
0 & \quad \text{otherwise},
\end{cases}\\
\psi^*_k(\maya)&=
\begin{cases}
(-1)^{i}(\ldots,k_{i},k,k_{i+1},\ldots) & \quad \text{if}\  k_i<k<k_{i+1} \ \text{for some i},\\
0 & \quad \text{otherwise}, 
\end{cases}
\end{align*}
and $c=1$. This representation is called the {\bf fermionic Fock space}.

\subsubsection{}\label{bosonfermioncorrespondence}

Let us identify $\C\Z\otimes\B$ and $\F$ by the bijection $F$ given in \ref{bosonfermionmap}.
The action of $\heisen$ on $\B$ is naturally extended to the action on $\C\Z\otimes\B$.
Consider generating functions
\[
\Gamma^+(z)=\exp\left(\,\sum_{m=1}^{\infty}\frac{z^m}{m}p(-m)\right),\quad \Gamma^-(z)=\exp\left(\,\sum_{m=1}^{\infty}\frac{z^{-m}}{m}p(m)\right),
\]
\[
\Psi(z)=\sum_{j=-\infty}^{\infty}z^{j-\h{1}}\psi_j,\quad \Psi^*(z)=\sum_{j=-\infty}^{\infty}z^{-j-\h{1}}\psi^*_j.
\]
Let $e^K$, $D$ denote the shift and degree operators on $\C\Z\otimes\B$ defined as follows : 
\[
e^K([c]\otimes b)=[c+1]\otimes b,\quad D([c]\otimes b)=c([c]\otimes b)\quad ([c]\otimes b\in\C\Z\otimes\B).
\]

\begin{thm}(See \cite{kac} Theorem14.10)
\begin{align*}
\Psi(z)&=\Gamma^+(z)\Gamma^-(z)^{-1}e^Kz^D,\\
\Psi^*(z)&=\Gamma^+(z)^{-1}\Gamma^-(z)e^{-K}z^{-D}.
\end{align*}
\end{thm}

\begin{rem}
Our notation is different from \cite{kac}'s. They do not use half integers. 
\end{rem}

\subsection{Frenkel-Kac construction}\label{subsection Frenkel-Kac construction}

\subsubsection{}\label{higherbosonrep}
Let $\heisen_{l-1}$ denote the $\C$-algebra generated by $p_i(m),p_i(-m)$ ($i\in\{1,\ldots,l-1\}$, $m>0$) and a central element $c$ with relations
%$k=\h{1},\ldots,l-\h{1}$, $m\in\Z_{>0}$ with relation
\[
[p_i(m),p_j(n)]=\delta_{0,n+m}a_{ij}c
\]
where $a_{ij}$ is the Cartan matrix of type $A_{l-1}$.
We say a representation of $\heisen_{l-1}$ is level-$1$ if $c$ acts by the identity map.

The affine Lie algebra of type $\hat{A}_{l-1}$ is $\hat{\mathfrak{sl}}_{l}=\mathfrak{sl}_{l}\otimes\C[x^{\pm 1}]\oplus\C c\oplus\C d$.
The subalgebra generated by $h_i\otimes x^m$ ($i\in\{1,\ldots,l-1\}$, $m\neq 0$) and $c$ is isomorphic to $\heisen_{l-1}$.

Let $Q$ be the root lattice of type $A_{l-1}$. 
We can regard $Q$ as the subset of $\ha$ by the isomorphism $\ha\simeq\ha^*$.
We can also identify $Q$ with $(\Z^{\tilde{I}})_0=\{(c_\h{1},\ldots,c_{l-\h{1}})\in\Z^{\tilde{I}}\mid \sum{c_k}=0\}$ so that simple roots 
$\alpha_i$ ($i\in\{1,\ldots,l-1\}$) correspond to ${\bf e}_{i-\h{1}}-{\bf e}_{i+\h{1}}$, where ${\bf e}_k$'s are the coordinate vector in $\Z^{\tilde{I}}$.
Then the set of positive roots is described as $P_+=\{{\bf e}_k-{\bf e}_{k'}\mid \h{1}\leq k<k'\leq l-\h{1}\}$.
For a positive root $\alpha={\bf e}_k-{\bf e}_{k'}$ 
let $e_{\alpha}$ (resp. $f_{\alpha}$) denote the matrix unit $E_{k+\h{1},k'+\h{1}}$ (resp. $E_{k'+\h{1},k+\h{1}}$) in $\mathfrak{sl}_{l}$.

\subsubsection{}\label{FrenkelKac}
Let $B$ be a level-$1$ $\heisen_{l-1}$-module. Assume that for any $b\in B$ there exists an integer $m(b)$ such that
\[
p_{k_1}(m_1)\cdots p_{k_a}(m_a)b=0\quad \text{if}\ m_i>0,\ \sum m_i>m(b).
\]
Set $F=\C Q \otimes B$. For $\alpha\in Q$, we define a generating function $X(\alpha,z)$ of operators on $F$ by
\[
X(\alpha,z)=\exp\left(\,\sum_{m=1}^{\infty}\frac{z^n}{n}\alpha(-n)\right)
\exp\left(-\sum_{m=1}^{\infty}\frac{z^{-n}}{n}\alpha(n)\right)
\exp\left(\log z\cdot\alpha(0)+\alpha\right)
\]
where $\exp(\log z\cdot\alpha(0)+\alpha)$ is the operator defined by
\[
\exp\left(\log z\cdot\alpha(0)+\alpha\right)([\beta] \otimes b)=z^{\h{1}(\alpha,\alpha)+(\alpha,\beta)}([\alpha+\beta] \otimes b)\quad ([\beta] \otimes b\in\C Q\otimes B).
\]
Let $X_m(\alpha)$ denote the operator given by 
\[
X(\alpha,z)=\sum_{m\in\Z}X_m(\alpha)z^m.
\]

We define a map $\varepsilon\colon Q\times Q \to\{\pm 1\}$ by 
\[
\varepsilon(\alpha_i,\alpha_j)=
\begin{cases}
-1 & (j=i,i+1)\\
1 & \text{otherwise}
\end{cases}
\]
and $\varepsilon(\alpha+\alpha',\beta)=\varepsilon(\alpha,\beta)\varepsilon(\alpha',\beta)$, 
$\varepsilon(\alpha,\beta+\beta')=\varepsilon(\alpha,\beta)\varepsilon(\alpha,\beta')$.

\begin{thm}(\cite{frenkel-kac})
The vector space $F=\C Q \otimes B$ has an $\hat{\mathfrak{sl}}_{l}$-module structure such that
\[
(h_i\otimes 1)\bigl([\beta] \otimes b\bigr)=\langle\alpha_i,\beta\rangle\bigl([\beta] \otimes b\bigr), \quad 
(h_i\otimes t^m)\bigl([\beta] \otimes b\bigr)=[\beta] \otimes p_i(m)b, \quad
\]
and for a positive root $\alpha$
\begin{align*}
(e_{\alpha}\otimes t^m)\bigl([\beta] \otimes b\bigr)&=\varepsilon(\alpha,\beta)X_m(\alpha)\bigl([\beta] \otimes b\bigr),\\
(f_{\alpha}\otimes t^m)\bigl([\beta] \otimes b\bigr)&=\varepsilon(\beta,\alpha)X_{-m}(-\alpha)\bigl([\beta] \otimes b\bigr) 
\end{align*}
and $c=1$, $d=0$.
\end{thm}

\begin{rem}
Identify $Q$ with $(\Z^{\tilde{I}})_0$ then we have 
\[
\varepsilon\left(\alpha_i,(c_\h{1},\ldots,c_{l-\h{1}})\right)=(-1)^{c_{i+\h{1}}},\quad
\varepsilon\left((c_\h{1},\ldots,c_{l-\h{1}}),\alpha_i)\right)=(-1)^{c_{i-\h{1}}}.
\]
\end{rem}

\subsection{Explicit formula for Frenkel-Kac construction}

\subsubsection{}
We define a level-$1$ $\heisen_{l-1}$-action on $\B^{\otimes {\tilde{I}}}$ by
\[
p_i(m)=p(m)_{i-\h{1}}-p(m)_{i+\h{1}}
\]
where $p(m)_k$ means the action of $p(m)$ on the $k$-th factor of $\B^{\otimes {\tilde{I}}}$.

We will give an explicit formula for the representation of $\hat{\mathfrak{sl}}_l$ 
obtained from this $\heisen_{l-1}$-action on $\B^{\otimes {\tilde{I}}}$ by Frenkel-Kac construction.

\subsubsection{}\label{formulaforFrenkelKac}

By the map $CQ$ defined in \ref{corequotient}, we can identify $\C Q\otimes\B^{\otimes \tilde{I}}$ and $\B$.

\begin{thm}
The action of $\hat{\mathfrak{sl}}_l$ on $\B$ 
obtained by Frenkel-Kac construction from the $\heisen_{l-1}$-module $\B^{\otimes \tilde{I}}$ is given as follows :
\begin{itemize}
\item $\hat{e}_i{b_\lambda}=(-1)^{v_i(\lambda)+v_{i+1}(\lambda)}\sum (-1)^{\eta^-(\lambda\backslash\mu,i,\lambda)}b_{\mu}$, where the summation runs over all $\mu$ obtained from $\lambda$ by removing a removable $i$-node,
\item $\hat{f}_i{b_\lambda}=(-1)^{v_{i-1}(\lambda)+v_{i}(\lambda)}\sum (-1)^{\eta^-(\mu\backslash\lambda,i,\lambda)}b_{\mu}$, where the summation runs over all $\mu$ obtained from $\lambda$ by adding an addable $i$-node, and
\item $\hat{h}_i{b_{\lambda}}=\left(|A_{\lambda,i}|-|R_{\lambda,i}|\right)b_{\lambda}$. 
\end{itemize}
\end{thm}
\begin{proof}
Here we only check the two actions of $\hat{e}_i$ coincide. 
Note that $\hat{e}_i=e_i\otimes 1\ (i=1,\ldots,l-1)$, $\hat{e}_0=\left(-\sum e_i\right)\otimes t$.
Since $p(m)_k$ and $p(m')_{k'}$ commute each other for any $m$ and $m'$ if $k\neq k'$, we have 
\begin{align*}
&\exp\left(\ \sum_{m=1}^{\infty}\frac{z^m}{m}\left(\p(-m)_{i-\h{1}}-\p(-m)_{i+\h{1}}\right)\right) 
\exp\left(-\sum_{m=1}^{\infty}\frac{z^{-m}}{m}\left(\p(m)_{i-\h{1}}-\p(m)_{i+\h{1}}\right)\right)\\
&=\Gamma^+_{i-\h{1}}(z)\circ\Gamma^+_{i+\h{1}}(z)^{-1}\circ\Gamma^-_{i-\h{1}}(z)^{-1}\circ\Gamma^-_{i+\h{1}}(z)\\
&=\Gamma^+_{i-\h{1}}(z)\circ\Gamma^-_{i-\h{1}}(z)^{-1}\circ\Gamma^+_{i+\h{1}}(z)^{-1}\circ\Gamma^-_{i+\h{1}}(z)\\
&=\left(\Psi_{i-\h{1}}(z)z^{-D_{i-\h{1}}}e^{-K_{i-\h{1}}}\right)\circ\left(\Psi^*_{i+\h{1}}(z)z^{D_{i+\h{1}}}e^{K_{i+\h{1}}}\right)
\quad\quad (\text{by Theorem \ref{bosonfermioncorrespondence}})\\
&=\Psi_{i-\h{1}}(z)\Psi^*_{i+\h{1}}(z)z^{-D_{i-\h{1}}}e^{-K_{i-\h{1}}}z^{D_{i+\h{1}}}e^{K_{i+\h{1}}}\\
&=\Psi_{i-\h{1}}(z)\Psi^*_{i+\h{1}}(z)z^{-1-D_{i-\h{1}}+D_{i+\h{1}}}e^{-K_{i-\h{1}}+K_{i+\h{1}}}.
\end{align*}
By definition we have
\[
\exp\left(\log z\cdot\alpha_i(0)+\alpha_i\right)=e^{K_{i-\h{1}}-K_{i+\h{1}}}z^{1+D_{i-\h{1}}-D_{i+\h{1}}}.
\]
Thus we get
\begin{align*}
X(\alpha_i,z)&=\Psi_{i-\h{1}}(z)\Psi^*_{i+\h{1}}(z)\,z\\
&=\sum_{m\in\Z}\left(\sum_{h\in\Z+\h{1}}\psi_{m+h,{i-\h{1}}}\psi^*_{h,{i+\h{1}}}\right)z^m,
%X_m(\alpha_i)=\sum_{j\in\Z+\h{1}}\psi_{m+j,{i-\h{1}}}\psi^*_{j,{i+\h{1}}},
\end{align*}
in particular, 
\begin{align*}
X_0(e_i)&=\sum_{h\in\Z+\h{1}}\psi_{h,{i-\h{1}}}\psi^*_{h,{i+\h{1}}},\\
X_1\left(-\sum e_i\right)&=\sum_{h\in\Z+\h{1}}\psi_{1+h,l-\h{1}}\psi^*_{h,\h{1}}.
\end{align*}
By \ref{eta} and \ref{Fermionicfockspace} we have  
\begin{align*}
&\psi_{h,i-\h{1}}\psi^*_{h,i+\h{1}}(b_\lambda)\\
=&\,\begin{cases}
\eta^-(X,i,\lambda)\cdot b_{\lambda\backslash X}&\left(\text{$\lambda$ has a removable node $X$ with content $l(h-\h{1})+i$}\right)\\
0&(\text{otherwise}),
\end{cases}
\end{align*}
and similar equations hold for $\psi_{1+h,l-\h{1}}\psi^*_{h,\h{1}}$.
Combine with Lemma \ref{corequotientrem}, then the claim follows.
\end{proof}

\section{Quiver varieties}\label{Quiver varieties}
{\ }

\vspace{-14pt}

In this section we study some properties of quiver varieties of type $\hat{A}$.

We have an $S^1$-action on the quiver varieties, so that the fixed points are isolated. 
The quiver varieties associated with different generic parameters are $S^1$-equivariantly diffeomorphic to each other,
and so there exists a natural bijection of fixed points.

For the construction of the representation of the affine Lie algebra, we use an "ordinary" parameter $\zeta_0$.
On the other hand, the Hilbert scheme of points on the minimal resolution corresponds to a parameter $\zeta_\infty$ 
in a chamber adjacent to the level-0 hyperplane.

The fixed points of the quiver varieties of type $\hat{A}$ associated with an "ordinary" parameter are parametrized by $\Pi$.
The fixed points of the Hilbert schemes of points on the minimal resolution of type ${A}$ are parametrized by $\Pi^l$.
The main results of this section is that the bijection described above is nothing but the operation $l$-core and $l$-quotient ({Theorem} \ref{thm1}).

\subsection{Quiver varieties}

\subsubsection{}
Let $(\I,H)$ be a quiver, namely $I$ is a set of vertices and $H$ is a set of oriented edges.
Assume we have a subset $\Omega\subset H$ 
such that $\Omega\cup\bar{\Omega}=H,\,\Omega\cap\bar{\Omega}=\varnothing$ where $\bar{\ }$ means reversing orientations of edges.
For 
$\dimv,\,\dimw\in\Z^{\I}_{\geq 0}$ and $\zeta=(\zeta_\C,\zeta_\R)\in \C^{\I}\oplus\R^{\I}$, 
we define the quiver variety $\M_{\zeta}(\dimv,\dimw)$ as follows (\cite{quiver1}). 
Here we assume $\dimv\in\Z^I_{>0}$, since otherwise we need a slight and non-essential modification.

Let $V$, $W$ be $\I$-graded vector spaces such that $\mathrm{dim}\,V_i=v_i$, $\mathrm{dim}\,W_i=w_i$.
Then we set
\[
\mathrm{M}(\dimv.\dimw)=
\left(\,\bigoplus_{h\in H}\Hom(V_{\mathrm{out}(h)},V_{\mathrm{in}(h)})\right)\oplus
\left(\,\bigoplus_{i\in I}\Hom(W_i,V_i)\oplus\Hom(V_i,W_i)\right), 
\]
where $h\in H$ is drawn from $\mathrm{out}(h)$ to $\mathrm{in}(h)$.
Note that $\GL_\dimv=\prod\GL(V_i)$ acts on $\mathrm{M}(\dimv,\dimw)$ by
\[
(g_i)\cdot\left(B_h,a_i,b_i\right)=\left(g_{\mathrm{in}(h)}B_hg_{\mathrm{out}(h)}^{-1},g_ia_i,b_ig_i^{-1}\right)
\]
for $g_i\in\GL(V_i)$, $B_h\in\Hom(V_{\mathrm{out}(h)},V_{\mathrm{in}(h)})$, $a_i\in\Hom(W_i,V_i)$ and $b_i\in\Hom(V_i,W_i)$.

The moment map $\mu_{\C}$ is given by 
\[
\mu_\C(B,a,b)=\left(\,\sum_{in(h)=i}\varepsilon(h)B_hB_{\bar{h}}+a_ib_i\right)
\in \bigoplus_{i\in I}\gl(V_i) = \gl_\dimv,
\]
where $\varepsilon(h)=1$ and $\varepsilon(\bar{h})=-1$ for $h\in\Omega$. 
Note that the center of $\gl_\dimv$ is canonically identified with $\C^{\I}$.

The quiver variety associated to $\dimv$, $\dimw$ and $\zeta$ is
\[
\M_\zeta(\dimv,\dimw)=
\left({\mu_\C}^{-1}(\zeta_\C)\right)^{ss}/\hspace{-3pt}/\,\GL_\dimv.
\]
By the index "ss", we mean the set of $\zeta_\R$-semistable objects in ${\mu_\C}^{-1}(\zeta_\C)$.
We omit the definition of $\zeta_\R$-semistability which is not required in the rest of this paper.

\subsubsection{}
Let $A=(a_{ij})_{i,j\in\I}$ be the adjacency matrix of the quiver $(I,H)$, that is to say 
$a_{ij}$ is the number of elements of $H$ which is drawn from $i$ to $j$.
Let $C=2\cdot\mathrm{id}-A$ be the Cartan matrix.
Then we consider the set of positive roots
\[
R_+ = \{\theta = (\theta_i)\in(Z_{\ge0})^{\I}\mid{}^t\theta C\theta\le 2\}\setminus\{0\},
\]
and for $\dimv\in\Z^{\I}$ we set 
\[
R_+(\dimv) = \{\theta\in R_+\mid\theta_i\le\dimv_i\ \forall i\in\I\}.
\]
The element $\zeta\in\C^I\oplus\R^I$ is called generic with respect to $\dimv$ 
if for any $\theta\in R_+(\dimv)$ 
\[
\zeta\notin D_{\theta}\otimes\R^3\subset\R^{\I}\otimes\R^3\simeq\C^{\I}\oplus\R^{\I}.
\]
where
\[
D_{\theta}=\left\{\zeta=(\zeta^i)\in\R^{\I}\ \big|\  \sum \zeta^i\theta_i=0\right\}\subset\R^{\I}.
\]

It is known that 
\begin{itemize}
\item $\M_{\zeta}(\dimv,\dimw)$ is smooth if $\zeta$ is generic with respect to $\dimv$.
Since $\R^{\I}\otimes\R^3\backslash \cup_\theta D_{\theta}\otimes\R^3$ is connected, $\M_{\zeta}(\dimv,\dimw)$ and $\M_{\zeta '}(\dimv,\dimw)$ are diffeomorphic for generic $\zeta$, $\zeta '$.
\item Let us fix $\zeta_\C\in\C^I$. Then we have a chamber structure on $\pi^{-1}(\zeta_{\C})$ defined by $D_{\theta}$'s, 
where $\pi$ is the projection $\C^I\oplus\R^I\to\C^I$.
If $\zeta$ and $\zeta '$ are in a same chamber, $\M_{\zeta}(\dimv,\dimw)$ and $\M_{\zeta '}(\dimv,\dimw)$ are 
isomorphic.
\item For generic $\zeta$, $\M_{\zeta}(\dimv,\dimw)$ is a fine moduli space of representations of the quiver.
Namely, on $\M_{\zeta}(\dimv,\dimw)$ there exists a universal family of representations of the quiver 
$(\mathfrak V, \mathfrak W, \mathfrak B, \mathfrak a, \mathfrak b)$.
\end{itemize}
Hereafter we sometimes identify a point of a quiver variety with the corresponding representation of the quiver.

\subsubsection{}\label{fixedpointsofquivervarieties}
We define an $S^1$-action on $\M_{\zeta}(\dimv,\dimw)$ by
\[
t*\left[\left(B_{\Omega},B_{\bar{\Omega}},a,b\right)\right]=\left[\left(tB_{\Omega},t^{-1}B_{\bar{\Omega}},a,b\right)\right].
\]
Note that this action may or may not be trivial.
For generic $\zeta$ and $\zeta '$, $\M_{\zeta}(\dimv,\dimw)$ and $\M_{\zeta '}(\dimv,\dimw)$ 
are $S^1$-equivariantly diffeomorphic.

Assume the $S^1$-fixed points are isolated for generic parameters.
Then for generic $\zeta$ and $\zeta '$ we have a canonical bijection between the $S^1$-fixed points 
of $\M_{\zeta}(\dimv,\dimw)$ and of $\M_{\zeta '}(\dimv,\dimw)$ induced by the $S^1$-equivariant diffeomorphism.

\subsubsection{}\label{fixedpointsofquivervarieties2}
For $\left[\left(B_{\Omega},B_{\bar{\Omega}},a,b\right)\right]\in\M_{\zeta}(\dimv,\dimw)^{S^1}$,
there exists a group homomorphism $\rho\colon S^1 \to\GL_\dimv$ such that
\[
t*\left(B_{\Omega},B_{\bar{\Omega}},a,b\right)=\rho(t)\left(B_{\Omega},B_{\bar{\Omega}},a,b\right)
\]
($\rho$ is uniquely determined because of the stability condition).
This induces an $S^1$-action on $V$.
The bijection described above preserve $V\in R(S^1)$ where $R(S^1)\simeq \Z[z^{\pm}]$ is the representation ring of $S^1$.

\subsection{Minimal resolution}
Let $(I,\Omega)$ be the quiver of type $\hat{A}_{l-1}$ with cyclic orientation. 
We set $\dimw_0=(1,0,\ldots,0)$, $\delta=(1,\ldots,1)\in\Z^{\I}$.

From now on we will consider the quiver varieties $\M_{(0,\zeta_\R)(\dimv,\dimw_0)}$ of type $\hat{A}_{l-1}$ only.
We write simply $\zeta$ and $\M_{\zeta}(\dimv)$ for $(0, \zeta_\R)$ and $\M_{(0,\zeta_\R)(\dimv,\dimw_0)}$.

\subsubsection{}
\begin{thm}(\cite{kronheimer-ALE})
Assume $\zeta$ is generic with respect to $\delta$, then
$\quiv{}{\delta}$ is isomorphic to the minimal resolution of $\C^2/\hspace{-3pt}/(\Z/l\Z)$.
\end{thm}
\begin{rem}
For any generic $\zeta$'s, $\quiv{}{\delta}$'s are isomorphic to each other, 
but the universal family depends on the chamber.
\end{rem}

\subsubsection{}\label{minimalresolution-rep}

We set
\[
C_0=\left\{\zeta\in\R^{\I}\mid \zeta^i<0\quad \text{for all $i\in I$}\right\}.
\]
For $\zeta_0\in C_0$, $\zeta_0$-semistability condition is equivalent to 
that there exists no proper subspace of $V$ which is $B$-invariant and contains $\mathrm{Im}\,a$.
\begin{prop}
For $\zeta_{0}\in C_0$, 
$\quiv{0}{\delta}$ has $l$ $S^1$-fixed points $P_{\frac{1}{2}},P_{\h{3}},\ldots,P_{l-\h{1}}$ 
and $P_k$ corresponds to the following representation :
\begin{align*}
B_{i,i+1}&=
\begin{cases}
1 & (i=0,\ldots ,k-\h{3})\\
0 & (i=k-\h{1},\ldots ,l-1),
\end{cases}
\\ 
B_{i,i-1}&=
\begin{cases}
0 & (i=1,\ldots ,k+\h{1})\\
1 & (i=k+\h{3},\ldots ,l),
\end{cases}
\end{align*}
\vspace{-10pt}
\[
a=1
,\ 
b=0.
\]
Here we identity the $\Hom$ spaces with $\C$ by certain bases of $1$-dimensional spaces $V_i$ and $W_i$ 's. 
\end{prop}
\begin{proof}
For $i,j\in\Z$ we define $B_{i,j}\in\mathrm{Hom}(V_i,V_j)$
\[
B_{i,j}=
\begin{cases}
B_{j-1,j}\cdots B_{i,i+1} & (i<j), \\
1 & (i=j), \\
B_{j+1,j}\cdots B_{i,i-1} & (i>j).
\end{cases}
\]
Assume that $B_{0,k-\h{1}}\neq 0$ for $k\in\tilde{I}$.
%$B_{k-\h{3},k-\h{1}}B_{k-\h{5},k-\h{3}}\cdots B_{0,1}\neq 0$. 
We can consider the ratio
$B_{l,k-\h{1}}\big/B_{0,k-\h{1}}\in\C$
of the two elements in $\mathrm{Hom}(V_0,V_{k-\h{1}})=\mathrm{Hom}(V_l,V_{k-\h{1}})\simeq\C$.
For a fixed point we have
\begin{align*}
B_{l,k-\h{1}}\Big/B_{0,k-\h{1}}&=t*\left(B_{l,k-\h{1}}\Big/B_{0,k-\h{1}}\right)\\
&=t^l\cdot B_{l,k-\h{1}}\Big/B_{0,k-\h{1}}\quad(t\in S^1),
\end{align*}
and so $B_{l,k-\h{1}}$ should be $0$.

By the parallel argument we also have $B_{0,l}=0$. 
So without loss of generality we may assume further $B_{0,k+\h{1}}=0$. Since $B_{0,k-\h{1}}\neq 0$ we have $B_{k-\h{1},k+\h{1}}=0$. 

The stability condition assures $B_{l,k+\h{1}}\neq 0$. Since $B_{l,k-\h{1}}=0$ we have $B_{k+\h{1},k-\h{1}}=0$. 

Vanishing of the value of the moment map assures  
$B_{i,i+1}=0$ ($i=k+\h{1},\ldots ,l-1$) and
$B_{i,i-1}=0$ ($i=1,\ldots ,k-\h{1}$).

By the $\GL_\dimv$-action we get the normalized forms as in the statement.
\end{proof}

\subsubsection{}\label{minimalresolution-coord}
\begin{prop}
We have a coordinate $(x,y)$ of $\quiv{0}{\delta}$ on a neighberhood of $P_k$ such that 
the representation corresponding to a point $(x,y)$ is given by
\begin{align*}
B_{i,i+1}&=
\begin{cases}
1 & (i=0,\ldots ,k-\h{3})\\
x & (i=k-\h{1})\\
xy & (i=k+\h{1},\ldots ,n-1),
\end{cases}
\\ 
B_{i,i-1}&=
\begin{cases}
xy & (i=1,\ldots ,k-\h{1})\\
y & (i=k+\h{1})\\
1 & (i=k+\h{3},\ldots ,n),
\end{cases}
\end{align*}
\vspace{-10pt}
\[ 
a=1
,\ 
b=0.
\]
In this coordinate, the $S^1$-action is described as $t*(x,y)=(t^lx.t^{-l}y)$.
\end{prop}
\begin{proof}
Note that $B_{i,i+1}$ ($i=0,\ldots ,k-\h{3}$) and $B_{i,i-1}$ ($i=k+\h{3},\ldots ,l$) do not vanish around $P_k$ 
and we can normalize them to $1$. 

Vanishing of the value of the moment map assures that the values of $B_{k\pm\h{1},k\mp\h{1}}$ determine the representation. 
We can see different values of $B_{k\pm\h{1},k\mp\h{1}}$ induce non-isomorphic representations. So the statement follows.
\end{proof}

\subsubsection{}\label{picofminimal}
\begin{exa}
(\,the case $l=4$\,)

The following figure exhibits the representations corresponding to the fixed point $P_{\h{3}}$ and 
a point $(x,y)$ in the local coordinate around $P_{\h{3}}$.

The three curves exhibit the exceptional curves in the minimal resolution, 
and the four dots exhibit the fixed points.

We have five dots in each box.
The top dot represents the basis of $W$ and the others represent the basis of $V$. 
The numbering of the vertices of the quiver is given counterclockwise. 

\begin{center}
%WinTpicVersion3.08
\unitlength 0.1in
\begin{picture}( 34.8000, 30.4000)(  8.0000,-37.0000)
% ELLIPSE 2 0 3 0
% 4 1400 1200 2000 1400 2000 1200 1200 1200
% 
\special{pn 8}%
\special{ar 1400 1200 600 200  3.1415927 6.2831853}%
% ELLIPSE 2 0 3 0
% 4 2200 1200 2800 1400 2800 1200 2000 1200
% 
\special{pn 8}%
\special{ar 2200 1200 600 200  3.1415927 6.2831853}%
% ELLIPSE 2 0 3 0
% 4 3000 1200 3600 1400 3600 1200 2800 1200
% 
\special{pn 8}%
\special{ar 3000 1200 600 200  3.1415927 6.2831853}%
% DOT 0 0 3 0
% 2 3400 1050 3400 1050
% 
\special{pn 20}%
\special{sh 1}%
\special{ar 3400 1050 10 10 0  6.28318530717959E+0000}%
\special{sh 1}%
\special{ar 3400 1050 10 10 0  6.28318530717959E+0000}%
% DOT 0 0 3 0
% 2 2600 1050 2600 1050
% 
\special{pn 20}%
\special{sh 1}%
\special{ar 2600 1050 10 10 0  6.28318530717959E+0000}%
\special{sh 1}%
\special{ar 2600 1050 10 10 0  6.28318530717959E+0000}%
% DOT 0 0 3 0
% 2 1800 1050 1800 1050
% 
\special{pn 20}%
\special{sh 1}%
\special{ar 1800 1050 10 10 0  6.28318530717959E+0000}%
\special{sh 1}%
\special{ar 1800 1050 10 10 0  6.28318530717959E+0000}%
% DOT 0 0 3 0
% 2 1000 1050 1000 1050
% 
\special{pn 20}%
\special{sh 1}%
\special{ar 1000 1050 10 10 0  6.28318530717959E+0000}%
\special{sh 1}%
\special{ar 1000 1050 10 10 0  6.28318530717959E+0000}%
% DOT 0 0 3 0
% 6 1800 2400 1200 3000 1800 3600 2400 3000 1800 2000 1800 2000
% 
\special{pn 20}%
\special{sh 1}%
\special{ar 1800 2400 10 10 0  6.28318530717959E+0000}%
\special{sh 1}%
\special{ar 1200 3000 10 10 0  6.28318530717959E+0000}%
\special{sh 1}%
\special{ar 1800 3600 10 10 0  6.28318530717959E+0000}%
\special{sh 1}%
\special{ar 2400 3000 10 10 0  6.28318530717959E+0000}%
\special{sh 1}%
\special{ar 1800 2000 10 10 0  6.28318530717959E+0000}%
\special{sh 1}%
\special{ar 1800 2000 10 10 0  6.28318530717959E+0000}%
% DOT 0 0 3 0
% 6 3600 2400 3000 3000 3600 3600 4200 3000 3600 2000 3600 2000
% 
\special{pn 20}%
\special{sh 1}%
\special{ar 3600 2400 10 10 0  6.28318530717959E+0000}%
\special{sh 1}%
\special{ar 3000 3000 10 10 0  6.28318530717959E+0000}%
\special{sh 1}%
\special{ar 3600 3600 10 10 0  6.28318530717959E+0000}%
\special{sh 1}%
\special{ar 4200 3000 10 10 0  6.28318530717959E+0000}%
\special{sh 1}%
\special{ar 3600 2000 10 10 0  6.28318530717959E+0000}%
\special{sh 1}%
\special{ar 3600 2000 10 10 0  6.28318530717959E+0000}%
% VECTOR 2 0 3 0
% 8 3000 2960 3560 2400 4200 2960 3640 2400 3640 3600 4200 3040 3000 3040 3560 3600
% 
\special{pn 8}%
\special{pa 3000 2960}%
\special{pa 3560 2400}%
\special{fp}%
\special{sh 1}%
\special{pa 3560 2400}%
\special{pa 3500 2434}%
\special{pa 3522 2438}%
\special{pa 3528 2462}%
\special{pa 3560 2400}%
\special{fp}%
\special{pa 4200 2960}%
\special{pa 3640 2400}%
\special{fp}%
\special{sh 1}%
\special{pa 3640 2400}%
\special{pa 3674 2462}%
\special{pa 3678 2438}%
\special{pa 3702 2434}%
\special{pa 3640 2400}%
\special{fp}%
\special{pa 3640 3600}%
\special{pa 4200 3040}%
\special{fp}%
\special{sh 1}%
\special{pa 4200 3040}%
\special{pa 4140 3074}%
\special{pa 4162 3078}%
\special{pa 4168 3102}%
\special{pa 4200 3040}%
\special{fp}%
\special{pa 3000 3040}%
\special{pa 3560 3600}%
\special{fp}%
\special{sh 1}%
\special{pa 3560 3600}%
\special{pa 3528 3540}%
\special{pa 3522 3562}%
\special{pa 3500 3568}%
\special{pa 3560 3600}%
\special{fp}%
% STR 2 0 3 0
% 3 1720 2100 1720 2200 5 0
% $1$
\put(17.2000,-22.0000){\makebox(0,0){$1$}}%
% STR 2 0 3 0
% 3 1440 2540 1440 2640 5 0
% $1$
\put(14.4000,-26.4000){\makebox(0,0){$1$}}%
% STR 2 0 3 0
% 3 2160 2540 2160 2640 5 0
% $1$
\put(21.6000,-26.4000){\makebox(0,0){$1$}}%
% STR 2 0 3 0
% 3 2160 3260 2160 3360 5 0
% $1$
\put(21.6000,-33.6000){\makebox(0,0){$1$}}%
% STR 2 0 3 0
% 3 3820 2680 3820 2780 5 0
% $1$
\put(38.2000,-27.8000){\makebox(0,0){$1$}}%
% STR 2 0 3 0
% 3 3380 2680 3380 2780 5 0
% $1$
\put(33.8000,-27.8000){\makebox(0,0){$1$}}%
% STR 2 0 3 0
% 3 3820 3120 3820 3220 5 0
% $1$
\put(38.2000,-32.2000){\makebox(0,0){$1$}}%
% STR 2 0 3 0
% 3 3380 3120 3380 3220 5 0
% $y$
\put(33.8000,-32.2000){\makebox(0,0){$y$}}%
% STR 2 0 3 0
% 3 3220 3280 3220 3380 5 0
% $x$
\put(32.2000,-33.8000){\makebox(0,0){$x$}}%
% STR 2 0 3 0
% 3 4020 3280 4020 3380 5 0
% $xy$
\put(40.2000,-33.8000){\makebox(0,0){$xy$}}%
% STR 2 0 3 0
% 3 4010 2520 4010 2620 5 0
% $xy$
\put(40.1000,-26.2000){\makebox(0,0){$xy$}}%
% STR 2 0 3 0
% 3 3190 2520 3190 2620 5 0
% $xy$
\put(31.9000,-26.2000){\makebox(0,0){$xy$}}%
% CIRCLE 2 2 3 0
% 4 1800 1060 1800 1460 1800 1460 1800 1460
% 
\special{pn 8}%
\special{ar 1800 1060 400 400  0.0000000 0.0300000}%
\special{ar 1800 1060 400 400  0.1200000 0.1500000}%
\special{ar 1800 1060 400 400  0.2400000 0.2700000}%
\special{ar 1800 1060 400 400  0.3600000 0.3900000}%
\special{ar 1800 1060 400 400  0.4800000 0.5100000}%
\special{ar 1800 1060 400 400  0.6000000 0.6300000}%
\special{ar 1800 1060 400 400  0.7200000 0.7500000}%
\special{ar 1800 1060 400 400  0.8400000 0.8700000}%
\special{ar 1800 1060 400 400  0.9600000 0.9900000}%
\special{ar 1800 1060 400 400  1.0800000 1.1100000}%
\special{ar 1800 1060 400 400  1.2000000 1.2300000}%
\special{ar 1800 1060 400 400  1.3200000 1.3500000}%
\special{ar 1800 1060 400 400  1.4400000 1.4700000}%
\special{ar 1800 1060 400 400  1.5600000 1.5900000}%
\special{ar 1800 1060 400 400  1.6800000 1.7100000}%
\special{ar 1800 1060 400 400  1.8000000 1.8300000}%
\special{ar 1800 1060 400 400  1.9200000 1.9500000}%
\special{ar 1800 1060 400 400  2.0400000 2.0700000}%
\special{ar 1800 1060 400 400  2.1600000 2.1900000}%
\special{ar 1800 1060 400 400  2.2800000 2.3100000}%
\special{ar 1800 1060 400 400  2.4000000 2.4300000}%
\special{ar 1800 1060 400 400  2.5200000 2.5500000}%
\special{ar 1800 1060 400 400  2.6400000 2.6700000}%
\special{ar 1800 1060 400 400  2.7600000 2.7900000}%
\special{ar 1800 1060 400 400  2.8800000 2.9100000}%
\special{ar 1800 1060 400 400  3.0000000 3.0300000}%
\special{ar 1800 1060 400 400  3.1200000 3.1500000}%
\special{ar 1800 1060 400 400  3.2400000 3.2700000}%
\special{ar 1800 1060 400 400  3.3600000 3.3900000}%
\special{ar 1800 1060 400 400  3.4800000 3.5100000}%
\special{ar 1800 1060 400 400  3.6000000 3.6300000}%
\special{ar 1800 1060 400 400  3.7200000 3.7500000}%
\special{ar 1800 1060 400 400  3.8400000 3.8700000}%
\special{ar 1800 1060 400 400  3.9600000 3.9900000}%
\special{ar 1800 1060 400 400  4.0800000 4.1100000}%
\special{ar 1800 1060 400 400  4.2000000 4.2300000}%
\special{ar 1800 1060 400 400  4.3200000 4.3500000}%
\special{ar 1800 1060 400 400  4.4400000 4.4700000}%
\special{ar 1800 1060 400 400  4.5600000 4.5900000}%
\special{ar 1800 1060 400 400  4.6800000 4.7100000}%
\special{ar 1800 1060 400 400  4.8000000 4.8300000}%
\special{ar 1800 1060 400 400  4.9200000 4.9500000}%
\special{ar 1800 1060 400 400  5.0400000 5.0700000}%
\special{ar 1800 1060 400 400  5.1600000 5.1900000}%
\special{ar 1800 1060 400 400  5.2800000 5.3100000}%
\special{ar 1800 1060 400 400  5.4000000 5.4300000}%
\special{ar 1800 1060 400 400  5.5200000 5.5500000}%
\special{ar 1800 1060 400 400  5.6400000 5.6700000}%
\special{ar 1800 1060 400 400  5.7600000 5.7900000}%
\special{ar 1800 1060 400 400  5.8800000 5.9100000}%
\special{ar 1800 1060 400 400  6.0000000 6.0300000}%
\special{ar 1800 1060 400 400  6.1200000 6.1500000}%
\special{ar 1800 1060 400 400  6.2400000 6.2700000}%
% LINE 2 1 3 0
% 2 1800 1110 1800 1900
% 
\special{pn 8}%
\special{pa 1800 1110}%
\special{pa 1800 1900}%
\special{da 0.070}%
% BOX 2 1 3 0
% 2 1080 1900 2480 3700
% 
\special{pn 8}%
\special{pa 1080 1900}%
\special{pa 2480 1900}%
\special{pa 2480 3700}%
\special{pa 1080 3700}%
\special{pa 1080 1900}%
\special{da 0.070}%
% VECTOR 2 0 3 0
% 2 3600 2440 3040 3000
% 
\special{pn 8}%
\special{pa 3600 2440}%
\special{pa 3040 3000}%
\special{fp}%
\special{sh 1}%
\special{pa 3040 3000}%
\special{pa 3102 2968}%
\special{pa 3078 2962}%
\special{pa 3074 2940}%
\special{pa 3040 3000}%
\special{fp}%
% VECTOR 2 0 3 0
% 2 3600 3560 3040 3000
% 
\special{pn 8}%
\special{pa 3600 3560}%
\special{pa 3040 3000}%
\special{fp}%
\special{sh 1}%
\special{pa 3040 3000}%
\special{pa 3074 3062}%
\special{pa 3078 3038}%
\special{pa 3102 3034}%
\special{pa 3040 3000}%
\special{fp}%
% VECTOR 2 0 3 0
% 2 4160 3000 3600 3560
% 
\special{pn 8}%
\special{pa 4160 3000}%
\special{pa 3600 3560}%
\special{fp}%
\special{sh 1}%
\special{pa 3600 3560}%
\special{pa 3662 3528}%
\special{pa 3638 3522}%
\special{pa 3634 3500}%
\special{pa 3600 3560}%
\special{fp}%
% VECTOR 2 0 3 0
% 2 3600 2440 4160 3000
% 
\special{pn 8}%
\special{pa 3600 2440}%
\special{pa 4160 3000}%
\special{fp}%
\special{sh 1}%
\special{pa 4160 3000}%
\special{pa 4128 2940}%
\special{pa 4122 2962}%
\special{pa 4100 2968}%
\special{pa 4160 3000}%
\special{fp}%
% VECTOR 2 0 3 0
% 2 3600 2040 3600 2360
% 
\special{pn 8}%
\special{pa 3600 2040}%
\special{pa 3600 2360}%
\special{fp}%
\special{sh 1}%
\special{pa 3600 2360}%
\special{pa 3620 2294}%
\special{pa 3600 2308}%
\special{pa 3580 2294}%
\special{pa 3600 2360}%
\special{fp}%
% VECTOR 2 0 3 0
% 2 1800 2040 1800 2360
% 
\special{pn 8}%
\special{pa 1800 2040}%
\special{pa 1800 2360}%
\special{fp}%
\special{sh 1}%
\special{pa 1800 2360}%
\special{pa 1820 2294}%
\special{pa 1800 2308}%
\special{pa 1780 2294}%
\special{pa 1800 2360}%
\special{fp}%
% VECTOR 2 0 3 0
% 2 1770 2430 1230 2970
% 
\special{pn 8}%
\special{pa 1770 2430}%
\special{pa 1230 2970}%
\special{fp}%
\special{sh 1}%
\special{pa 1230 2970}%
\special{pa 1292 2938}%
\special{pa 1268 2932}%
\special{pa 1264 2910}%
\special{pa 1230 2970}%
\special{fp}%
% VECTOR 2 0 3 0
% 2 2370 3030 1830 3570
% 
\special{pn 8}%
\special{pa 2370 3030}%
\special{pa 1830 3570}%
\special{fp}%
\special{sh 1}%
\special{pa 1830 3570}%
\special{pa 1892 3538}%
\special{pa 1868 3532}%
\special{pa 1864 3510}%
\special{pa 1830 3570}%
\special{fp}%
% VECTOR 2 0 3 0
% 2 1830 2430 2370 2970
% 
\special{pn 8}%
\special{pa 1830 2430}%
\special{pa 2370 2970}%
\special{fp}%
\special{sh 1}%
\special{pa 2370 2970}%
\special{pa 2338 2910}%
\special{pa 2332 2932}%
\special{pa 2310 2938}%
\special{pa 2370 2970}%
\special{fp}%
% BOX 2 1 3 0
% 2 2880 1900 4280 3700
% 
\special{pn 8}%
\special{pa 2880 1900}%
\special{pa 4280 1900}%
\special{pa 4280 3700}%
\special{pa 2880 3700}%
\special{pa 2880 1900}%
\special{da 0.070}%
% LINE 2 1 3 0
% 2 3600 1900 2000 1300
% 
\special{pn 8}%
\special{pa 3600 1900}%
\special{pa 2000 1300}%
\special{da 0.070}%
% STR 2 0 3 0
% 3 3520 2100 3520 2200 5 0
% $1$
\put(35.2000,-22.0000){\makebox(0,0){$1$}}%
\end{picture}%
\end{center}
\end{exa}

\subsection{$\Z/l\Z$-equivariant Hilbert scheme ($\zeta_0$-case)}

\subsubsection{}\label{equivHilb}

Let $(\C^2)^{[n]}$ denote the Hilbert scheme of $n$ points on $\C^2$ : 
\[
\hilb=\{J\underset{\text{ideal}}{\subset}\C[z_1,z_2]\mid \mathrm{dim}\,\C[z_1,z_2]/J=n\}.
\] 

The cyclic group $\Z/l\Z$ acts on $\C^2$ by $r\cdot(z_1,z_2)=(rz_1,^{-1}z_2)$ ($r\in\Z/l\Z$). 
This action induces an action of $\Z/l\Z$ on $(\C^2)^{[n]}$.
Let $\left((\C^2)^{[n]}\right)^{\Z/l\Z}$ be the set of fixed points.
For $J\in\left((\C^2)^{[n]}\right)^{\Z/l\Z}$ , 
$\C[z_1,z_2]/J$ has a canonical $\Z/l\Z$-module structure. 
Let $\left((\C^2)^{[n]}\right)^{\Z/l\Z,\,\dimv}\subset\left((\C^2)^{[n]}\right)^{\Z/l\Z}$ 
denote the set of points such that the corresponding $\Z/l\Z$-module is isomorphic to $\oplus_i(\C_{(i)})^{\oplus \dimv_i}$,
where $\C_{(i)}$ is the 1-dimensional representation of $\Z/l\Z$ with weight $i$.
\begin{thm}(\cite{nakajima-lec-note} Theorem 4.4)
\begin{enumerate}
\item[(i)]
For $\zeta_0\in C_0$, $\quiv{0}{\dimv}$ is isomorphic to $\left((\C^2)^{[n]}\right)^{\Z/l\Z,\dimv}$. 
\item[(ii)]
The representation $(V(J),W(J),B(J),a(J),b(J))$ of the quiver corresponding to $J\in\left((\C^2)^{[n]}\right)^{\Z/l\Z}$ is given as follows : 
\[
V(J)_i=\mathrm{Hom}_{\Z/l\Z}(\C_{(i)},\C[z_1,z_2]/J),\quad\quad\quad
W(J)_0=\C
\]
and
\[
\begin{array}{c}
B(J)_{i,i+1}=z_1|_{V_i}\colon V_i\to V_{i+1},\quad\quad
B(J)_{i,i-1}=z_2|_{V_i}\colon V_i\to V_{i-1},\smallskip\\
a(J)(1)=[1]\in \C[z_1,z_2]/J,\quad\quad
b(J)=0,
\end{array}
\]
where $z_1$ and $z_2$ represents the multiple operators on $\C[z_1,z_2]/J$.
\end{enumerate}
\end{thm}

\subsubsection{}\label{paraoffixpt0}
Let us consider the $S^1$-action on $\C^2$ given by $t\cdot(z_1,z_2)=(tz_1,t^{-1}z_2)$. 
This induces a $S^1$-action on $(\C^2)^{[n]}$.
The embedding $\quiv{0}{\dimv}\simeq\left((\C^2)^{[n]}\right)^{\Z/l\Z,\,\dimv}\subset(\C^2)^{[n]}$ is $S^1$-equivariant. 
So we get $\quiv{0}{\dimv}^{S^1}=\quiv{0}{\dimv}\cup\left((\C^2)^{[n]}\right)^{S^1}$.

The $S^1$-fixed points of $(\C^2)^{[n]}$ are parametrized by Young diagrams (\cite{nakajima-lec-note}).
For $\lambda\in\Pi$ the corresponding ideal $J_\lambda\in \left(\C^2\right)^{\left[|\lambda|\right]}$ is the ideal generated by $\{{z_1}^a{z_2}^b\mid (a,b)\notin\lambda\}$.
Then $\{[{z_1}^a{z_2}^b]\in\C[z_1,z_2]/J_\lambda\mid (a,b)\in\lambda\}$ forms a basis of $\C[z_1,z_2]/J_\lambda$. 
Since $r\cdot[{z_1}^a{z_2}^b]=r^{a-b}[{z_1}^a{z_2}^b]$ for $r\in\Z/l\Z$ we have $\quiv{0}{\dimv}^{S^1}=\{J_\lambda\mid\dimv(\lambda)=\dimv\}$.

\subsubsection{}\label{repassfixpt0}
\begin{prop}
The representation $(V(J_\lambda),W(J_\lambda),B(J_\lambda),a(J_\lambda),b(J_\lambda))$ of the quiver corresponding to a fixed point $J_\lambda$ is described as follows :  
\[
V(J_\lambda)_i=\bigoplus_{a-b\equiv i,(a,b)\in\lambda}\C\,v_{(a,b)}, \quad W(J_\lambda)_0=\C\,w,
\]
and
\begin{align*}
B(J_\lambda)_{i.i+1}(v_{(a,b)})&=
\begin{cases}
v_{(a+1,b)} & \text{if}\ (a+1,b)\in\lambda,\\
0 & \text{if}\ (a+1,b)\notin\lambda,
\end{cases}\\
B(J_\lambda)_{i.i-1}(v_{(a,b)})&=
\begin{cases}
v_{(a,b+1)} & \text{if}\ (a,b+1)\in\lambda,\\
0 & \text{if}\ (a,b+1)\notin\lambda,
\end{cases}\\
a(J_\lambda)(w)&=v_{(0,0)},\quad b(J_\lambda)=0.
\end{align*}
\end{prop}
\begin{proof}
Apply Theorem \ref{equivHilb} (ii) for $\C[x,y]/J_\lambda$.
\end{proof}
\begin{cor}
\[V(J_\lambda)=f_\lambda(z)\in R(S^1).\]
\end{cor}
\begin{proof}
We can see $\rho(t)(v_{(a,b)})=t^{a-b}v_{(a,b)}$. So the claim follows.
\end{proof}

\subsubsection{}\label{tangentspace0}
The character of the tangent space of $(\C^2)^{[n]}$ at fixed point $J_\lambda$ is 
\[
T_{J_\lambda}(\C^2)^{[n]}=\sum_{(a,b)\in\lambda}t^{h(a,b)}+t^{-h(a,b)}
\] 
where $h(a,b)$ is the hook length of the hook associated with a node $(a,b)$ (\cite{nakajima-lec-note}).
Further we can see
\[
T_{J_\lambda}\quiv{0}{\dimv}=\sum_{(a,b)}t^{h(a,b)}+t^{-h(a,b)}
\] 
where the summation runs over all $(a,b)\in\lambda$ such that $h(a,b)\equiv 0\ (\mathrm{mod}\ l)$

\begin{rem}
In term of Maya diagrams, such a hook corresponds to 
\[
(k,k+nl)\in(\Z+1/2)^2\quad\text{such that}\ \ n>0,\ \maya(k)=1,\ \maya(k+nl)=-1 
\]
and its hook length is $nl$.
\end{rem}

\subsection{Hilbert scheme of the minimal resolution ($\zeta_\infty$-case)}

\subsubsection{}\label{Hilbofminimalresolution}
The hyperplane $D_\delta=\{\zeta\in \R^I\mid \sum\zeta_i=0\}$ of $\R^I$ is called the level $0$ hyperplane. 
We have the chamber structure on this hyperplane defined by $D_\alpha$'s, where $\alpha$ is a root of finite root system, and
\[
C=\{\zeta\in D_\delta\mid \zeta_i>0\ (i=1,\ldots,l-1)\}\subset D_\delta
\]
is one of the chambers.

For $\dimv\in\Z^l$ let $C_\infty(\dimv)$ denote the unique chamber which is contained in $\{\zeta\in\R^I\mid \sum\zeta_i>0\}$ and 
has $C$ as its face.

Let $M$ denote the minimal resolution of $\C^2/\hspace{-3pt}/(\Z/l\Z)$. 
For $\zeta_{\infty}\in C_{\infty}(\dimv)$ the quiver variety $\quiv{\infty}{\dimv}$ is a certain moduli space of torsion free sheaves on $M$. 
In next subsection we will review this result very briefly. 
The reader can refer \cite{nakajima-moduli-on-ALE} for detail.

\subsubsection{}\label{}
Let us consider the action of $\Z/l\Z$ on $\CP^2$ given by $r[x:y:z]=[rx:r^{-1}y:z]$.
Let $\tilde{M}$ denote the orbifold which is a compactification of $M$ given by resolving the singular point $\left[[0:0:1]\right]$ of $\CP^2/\hspace{-3pt}/(\Z/l\Z)$.
Note that $\tilde{M}$ has the natural $S^1$-action compatible with the one on $M$.
We set $l_\infty=\tilde{M}\backslash M$.

Let $\V=\oplus_i\V_i$ be the universal bundle on $M$ and $\tV=\oplus_i\tV_i$ be its extension to $\tM$.
Note that $\tV$ has the unique $S^1$-equivariant structure such that the restriction of the action to $l_\infty$ is trivial.

\begin{thm}
For $\zeta_{\infty}\in C_{\infty}(\dimv)$ the quiver variety $\quiv{\infty}{\dimv}$ is the fine moduli space of rank 1 torsion free sheaves 
$E$ on $\tM$ such that $E|_{l_\infty}$ is trivial and $c_1(E)$ and $\mathrm{ch}_2(E)$ take values given in (1.8) of \cite{nakajima-moduli-on-ALE}.   
\end{thm}

\begin{rem}
(i) For a representation $(V,W,B_\Omega,B_{\bar{\Omega}},a,b)$, corresponding torsion free sheaf $E(V,W,B_\Omega,B_{\bar{\Omega}},a,b)$ 
is given as the cohomology of the following complex of sheaves :
\[
C^\bullet (V,W,B_\Omega,B_{\bar{\Omega}},a,b)\colon\begin{array}{ccccc}
&& V\otimes\tV &&\\
&& \oplus &&\\
V\otimes\tV(-l_\infty) & 
\overset{\alpha}{\longrightarrow} & 
V\otimes\tV &
\overset{\beta}{\longrightarrow} &
V\otimes\tV(l_\infty)\\
&& \oplus &&\\
&& W\otimes\tV &&
\end{array}
\]
where
\[
\alpha=\left(
\begin{array}{c}
B_{\Omega}\otimes \tilde{1}-1\otimes \tilde{\mathfrak{B}}_{{\Omega}}\\
B_{\bar{\Omega}}\otimes \tilde{1}-1\otimes \tilde{\mathfrak{B}}_{\bar{\Omega}}\\
a\otimes \tilde{1}
\end{array}
\right),\quad
\beta=(-B_{\bar{\Omega}}\otimes \tilde{1}+1\otimes \tilde{\mathfrak{B}}_{\bar{\Omega}},B_{{\Omega}}\otimes \tilde{1}-1\otimes \tilde{\mathfrak{B}}_{\bar{\Omega}},b\otimes \tilde{1}).
\]
Here $\tilde{\mathfrak{B}}$ is the extension of universal family $\mathfrak{B}\colon \V\to\V$ of maps to $\tilde{M}$ 
and $\tilde{1}$ means the natural morphism $\tV(-l_\infty)\to\tV$ and $\tV\to\tV(l_\infty)$. 

\noindent (ii) For a torsion free sheaf $E$, the representation spaces $V(E)$ and $W(E)$ are $H^1(E\otimes \tV)$ and the fiber of $E$ on a point in $l_\infty$ respectively. 
Although we omit the constructions of the maps, we mention that they are functorial.  
\end{rem}

\subsubsection{}\label{qtoh}
For $E\in\quiv{\infty}{\dimv}$ its double dual $E^{\vee\vee}$ is a line bundle on $\tM$ such that its restriction to $l_\infty$ is trivial and 
$c_1(E^{\vee\vee})=c_1(E)$. Note that $E^{\vee\vee}$ is determined uniquely by these conditions. In fact, since  
$c_1(E)=\sum_{i\neq 0}\mathrm{u}_ic_1(\tV_i)$ where $\mathbf{u}=-C\dimv$ ((1.8) of \cite{nakajima-moduli-on-ALE})
and $\{c_1(\tV_i)\}$ is the dual basis of $\{[C_i]\}$ (\cite{kronheimer-nakajima} Proposition 2.2), we can check 
$E^{\vee\vee}=\mathcal{O}(\sum_{i=1}^{l-1}(v_i-v_0)C_i)$.

The quotient sheaf $E^{\vee\vee}/E$ is supported at finitely many points on $M$ and its length equals to $n=\mathrm{ch}_2(E^{\vee\vee})-\mathrm{ch}_2(E)$, 
which we can calculate from (1.8) of \cite{nakajima-moduli-on-ALE}. 
The map $E\mapsto E^{\vee\vee}/E$ induces the isomorphism between $\quiv{\infty}{\dimv}$ and the Hilbert scheme $M^{[n]}$ of $n$ points on $M$.
Since $\mathrm{dim}\,\quiv{\infty}{\dimv}=\dimv C\hspace{0.7pt}{}^t\hspace{-0.7pt}\dimv+2v_0$ we have $n=\dimv C\hspace{0.7pt}{}^t\hspace{-0.7pt}\dimv/2+v_0$.

For $c=(c_\h{1},\ldots,c_{l-\h{1}})\in Q$ and $n\in\Z_{\geq 0}$ let $\dimv\in\Z^I$ be the unique elements such that $c_k=v_{k-\h{1}}-v_{k+\h{1}}$ and $n=\dimv C\hspace{0.7pt}{}^t\hspace{-0.7pt}\dimv/2+v_0$.
Let $\phi_{c,n}=\phi_\dimv$ denote the isomorphism between $M^{[n]}$ and $\quiv{\infty}{\dimv}$.

\subsubsection{}\label{fixptinf}
Let $\mathrm{Sym}^n(M)$ be the $n$-th symmetric product of $M$ and $\pi\colon M^{[n]}\to\mathrm{Sym}^n(M)$ be the the Hilbert-Chow morphism.

The $S^1$-action on $M$ induces $S^1$-actions on $M^{[n]}$ and $\mathrm{Sym}^n(M)$ so that $\pi$ is $S^1$-equivariant.
We can see
\[
\mathrm{Sym}^n(M)^{S^1}=\left\{\sum n_k[P_k]\ \bigg|\  \sum n_k=n\right\}.
\]
There exists a neighborhood of $\sum n_k[P_k]$ which is isomorphic to some open set in $\prod\mathrm{Sym}^{n_k}(\C^2)$.
The inverse image for $\pi$ of this open set is isomorphic to some open set in $\prod\left(\C^2\right)^{[n_k]}$.
So $\left(M^{[n]}\right)^{S^1}$ is parametrized by $l$-tuple of Young diagrams $\blam=(\lambda_{\h{1}},\ldots\lambda_{l-\h{1}})\in\Pi^{\tilde{I}}$ 
such that $\sum|\lambda_k|=n$.
Let $J_{\blam}$ denote the corresponding fixed point.

\subsubsection{}
\begin{lem}
The $S^1$-action on the moduli space $\quiv{\infty}{\dimv}$ induced by the $S^1$-action on $\tM$ coincides with the one given in \ref{fixedpointsofquivervarieties}.  
\end{lem}
\begin{proof}
The $S^1$-equivariant structure of $\tV$ induces the isomorphism $\mathrm{Hom}(\tV,\tV)=\mathrm{Hom}(t^*\tV,t^*\tV)$ for $t\in S^1$.
Under this isomorphism we have $t^*\mathfrak{B}_\Omega=t^{-1}\mathfrak{B}_\Omega$ and $t^*\mathfrak{B}_{\bar{\Omega}}=t\mathfrak{B}_{\bar{\Omega}}$.
Thus the complex $C^\bullet(V,W,tB_\Omega,t^{-1}B_{\bar{\Omega}},a,b)$ 
is isomorphic to the complex $t^*C^\bullet(V,W,B_\Omega,B_{\bar{\Omega}},a,b)$.
Take the cohomology of the complex, then the claim follows.
\end{proof}
The isomorphism $\phi_{c,n}$ given in \ref{qtoh} is $S^1$-equivariant.
For $c\in Q$ and $\blam\in\Pi^{\tilde{I}}$, we set $E_{c,\blam}=\phi_{c,n}(J_{\blam})$. This gives a bijection between $Q\times \Pi^{\tilde{I}}$ and $\coprod_\dimv\quiv{\infty}{\dimv}^{S^1}$.

\subsubsection{}
For $E\in\quiv{\infty}{\dimv}^{S^1}$, $E^{\vee\vee}=\mathcal{O}(\sum_{i=1}^{l-1}(v_i-v_0)C_i)$ has the unique $S^1$-equivariant structure such that 
the restriction of the action on $l_\infty$ is trivial, and so is $E$. 

\begin{lem}
The $S^1$-action on $V(E)=H^1(E\otimes \tV)$ induced by the $S^1$-equivariant structures on $\tV$ and $E$ 
coincides with the one given in \ref{fixedpointsofquivervarieties2}.
\end{lem}
\begin{proof}
The $S^1$-action on $W(E)$ induced by the $S^1$-equivariant structures on $\tV$ and $E$ is trivial. 
Since every construction is functorial with respect to $E$, the homomorphism $S^1\to\mathrm{GL}(V(E))$ satisfies the condition of $\rho$ described in \ref{fixedpointsofquivervarieties2}.
\end{proof}

\subsubsection{}\label{fixptinfprop}
\begin{prop}
For $c\in Q$ we set $\lambda(c)=CQ^{-1}(c,(\emptyset,\ldots,\emptyset))$. Then we have
\[
H^1(E_{c,\blam}\otimes \tV)=f_{\lambda(c)}(z)+\sum_kz^{lc_k}f_{\lambda_k}(z^l)\left(z^{k-l+\h{1}}+\cdots+z^{k-\h{1}}\right)\in R(S^1).
\]
\end{prop}
\begin{proof}
The quiver variety $\quiv{\infty}{\dimv(\lambda(c))}$ is one point and the point corresponds to $E^{\vee\vee}$. 
So we have $H^1(E^{\vee\vee}\otimes \tV)=f_{\lambda(c)}(z)\in R(S^1)$.

Since $H^0(\tM,E^{\vee\vee}\otimes \tV)=0$ (\cite{nakajima-moduli-on-ALE} 5(ii)) we have the following exact sequence of $S^1$-module :
\[
0\to H^0(\tM,E^{\vee\vee}/E\otimes \tV)\to H^1(E\otimes \tV)\to H^1(E^{\vee\vee}\otimes \tV)\to 0.
\]
Let $(E^{\vee\vee}/E)_k$ denote the direct summand of $E^{\vee\vee}/E$ supported on $P_k$. 
On the neighborhood of $P_k$ given in \ref{minimalresolution-coord} we have
\begin{align*}
(E^{\vee\vee}/E)_k\otimes \tV
&=\mathcal{O}/J_{\lambda_k}\otimes E^{\vee\vee}\otimes \tV\\
%&\simeq \C[x,y]/J_{q_k(\lambda)}\otimes_{\C[x,y]}\C[x,y]\otimes_{\C[x,y]}\C[x,y]^{\oplus I}.
&\simeq \C[x,y]/J_{\lambda_k}\otimes E^{\vee\vee}_{\ \ P_k}\otimes \tV_{P_k}.
\end{align*}
Recall that the $S^1$-action on this neighborhood is given by $t*(x,y)=(t^lx,t^{-l}y)$. 
We have $\C[x,y]/J_{\lambda_k}=f_{\lambda_k}(z^l)$ as in \ref{repassfixpt0}.
Since $E^{\vee\vee}=\mathcal{O}((v_{k-\h{1}}-v_0)(\text{x-axis})+(v_{k+\h{1}}-v_0)(\text{y-axis}))$ on the neighborhood of $P_k$,
the weight of $E^{\vee\vee}_{\ \ P_k}$ is $l(v_{k-\h{1}}-v_{k+\h{1}})=lc_k$.
Using the description in \ref{minimalresolution-rep} we can see $\tV_{P_k}=z^{k-l+\h{1}}+\cdots+z^{k-\h{1}}$.
Then the claim follows.
\end{proof}

\subsection{Correspondence of fixed points}\label{thm1}
\begin{thm}
The following diagram is commutative :
\[
\begin{CD}
\coprod_{\dimv}\quiv{0}{\dimv}^{S^1} @>\ref{paraoffixpt0}>> \Pi\\
@V\ref{fixedpointsofquivervarieties}VV @VV\ref{corequotient}V\\
\coprod_{\dimv}\quiv{\infty}{\dimv}^{S^1} @>\ref{fixptinf}>> \ \ Q\times \Pi^{\tilde{I}}.
\end{CD}
\]
\end{thm}
\begin{proof}
Note that the map $\Pi\rightarrow\Z[z^\pm]$ given by $\lambda\mapsto f_\lambda(z)$
is injective. 
By \ref{fixedpointsofquivervarieties}, 
it is enough to check $V(J_\lambda)=V(E_{c(\lambda),q(\lambda)})\in R(S^1)$.
This follows from Proposition \ref{cqprop} and Proposition \ref{fixptinfprop}.
\end{proof}

\section{Representations on equivariant cohomologies}\label{Geometric representation assosiated with quiver varieties}
{\ }

\vspace{-14pt}

In this section we study the representations of the affine Lie algebra 
and the Heisenberg algebra on the middle degree $S^1$-equivariant cohomology groups of the quiver varieties.

First we see that the middle degree $S^1$-equivariant cohomology groups of quiver varieties 
has bases indexed by the $S^1$-fixed points (Proposition \ref{basis}).

The affine Lie algebra $\asl$ acts on $\oplus_{\dimv}H^{\mathrm{mid}}_{S^1}(\quiv{0}{\dimv})$ 
and the Heisenberg algebra $\heisen_{l-1}$ acts on $\oplus_{n}H^{\mathrm{mid}}_{S^1}(\quiv{\infty}{n\delta})$.
The main purpose of this section is to describe these actions with respect to above bases.

For the affine Lie algebra the argument works parallel with the one in \cite{vv-k}, 
dealing with the equivariant K-groups (Proposition \ref{calofVV}).
For the Heisenberg algebra we can find a formula in \cite{qin-wang} (Proposition \ref{thmofQW}). 
%This is a direct consequence of the results of \cite{Jack} or \cite{V}.    

\subsection{Equivariant cohomology groups}

We review some general results about $S^1$-equivariant cohomology groups. 
For more details, the reader can refer to \cite{audin} for example.

\subsubsection{}
We take $\C$ as the coefficient ring of cohomology groups.

Let $ES^1\to BS^1$ be the universal $S^1$-bundle which is given as 
the inductive limit of the Hopf fibration $S^{2n+1}\to \C P^{n}$. 
Note that the cohomology ring of $BS^1\simeq \C P^{\infty}$ is the 
polynomial ring $\C[t]$ with a generator $t\in H^2(BS^1)$.

For a topological space $X$ with an $S^1$-action, 
we define the $S^1$-equivariant cohomology group of $X$ by 
\[
H^*_{S^1}(X)=H^*(ES^1\times_{S^1}X).
\]

\subsubsection{}
Let $X$, $Y$ be $S^1$-equivariant topological spaces and $f\colon X\to Y$ be an $S^1$-equivariant map. 
Then we can define the following operators
\begin{align*}
\cup &\colon H^*_{S^1}(X)\otimes_\C H^*_{S^1}(X)\longrightarrow H^*_{S^1}(X)\quad (\text{cup product}),\\
f^*&\colon H^*_{S^1}(Y)\longrightarrow H^*_{S^1}(X)\quad (\text{pullback}).
\end{align*}
Note that pullbacks preserve cup products.

Let $p\colon X\to\{pt\}$. 
We have the action of $\C[t]=H^*_{S^1}(\{pt\})$ on $H^*_{S^1}(X)$ induced by $p_*$ and $\cup$.
We can see cup products and pullbacks commute with the $\C[t]$-actions.

\subsubsection{}

%Let $X$, $Y$ be smooth $S^1$-manifolds and $f\colon X\to Y$ be an $S^1$-invariant smooth map.

Let $Y$ a smooth $S^1$-manifold and $X$ be an $S^1$-invariant codimension $d$ smooth submanifold of $Y$.
Let $i\colon X\to Y$ denote the embedding. 
Then we can define a map 
\[
i_*\colon H^*_{S^1}(X)\to H^{*+d}_{S^1}(Y).
\]
(see \cite{audin} VI.4.c.). We set 
\[
[X]=i_*(1_X)\in H^{d}_{S^1}(Y)
\]
where $1_X\in H^{0}_{S^1}(X)$.

Let $W$ be a smooth $S^1$-manifold and  
$\pi\colon Z\to W$ be an $S^1$-equivariant fibre bundle whose fibre is a $d$-dimensional compact smooth manifold. 
Then we can define a map 
\[
\pi_*\colon H^*_{S^1}(Z)\to H^{*-d}_{S^1}(W).
\]
(see \cite{audin} VI.4.c.). 
Even in the case the fibre is not compact, if an element $\alpha\in H^*_{S^1}(Z)$ has the compact support on each fibre 
then we can define the element $\pi_*(\alpha)\in H^{*-d}_{S^1}(W)$. 
 
\begin{rem}
Roughly speaking, $i_*$ is induced by the Thom isomorphism with respect to the normal bundle on $X$ to $Y$ 
and $\pi_*$ is induced by the integration along fibres.
\end{rem}

\subsubsection{}\label{lemma for ec}
For an $S^1$-equivariant vector bundle $Z$ on an $S^1$-manifold $X$, we define the $S^1$-equivariant Euler class by
\[
\eu_{S_1}(Z)=\eu(ES_1\times_{S^1}Z)\in H^*(ES_1\times_{S^1} X)=H^*_{S_1}(X)
\]
where $\eu(\cdot)$ represents the Euler class of a vector bundle. 
To be precise, we define as the limit of $\eu(S^{2n+1}\times_{S^1}Z)$.
\begin{lem}
Let the all manifolds and morphisms below be smooth.
\begin{enumerate}
\item [(i)]  Let $i\colon X\to Y$ be an $S^1$-equivariant embedding and $\nu$ denote the normal bundle on $X$ to $Y$.
Then we have
\[
i^*i_*\alpha=\alpha\cup \eu_{S_1}(\nu)\quad (\alpha\in H^*_{S^1}(X)).
\]
\item [(ii)] (projection formula) 

Let $i\colon X\to Y$ be an $S^1$-equivariant embedding. Then we have
\[
i_*(\alpha\cup i^*\beta)=i_*\alpha\cup \beta\quad (\alpha\in H^*_{S^1}(X),\, \beta\in H^*_{S^1}(Y)).
\]
\item [(iii)] 
Let $i\colon X\to Y$ be an $S^1$-equivariant embedding and $Z$ be an $S^1$-manifold.

If $\alpha\in H^*_{S^1}(X\times Z)$ has the compact support on each fibre of $p_X\colon X\times Z\to X$, 
then $(i\times \mathrm{id})_*\alpha\in H^*_{S^1}(Y\times Z)$ also has the compact support on each fibre of $p_Y\colon Y\times Z\to Y$ and
\[
i_*({p_X}_*(\alpha))={p_Y}_*((i\times \mathrm{id})_*(\alpha)). 
\]
%\[
%\begin{CD}
%X\times Z   @>{i\times \mathrm{id}}>>   Y\times Z\\
%@V{p_X}VV                               @VV{p_Y}V\\
%X           @>{i}>>   Y   
%\end{CD}
%\]
\item [(iv)] (Thom class is Poincare dual of the zero section)

Let $\pi\colon Z\to W$ be an $S^1$-equivariant vector bundle and $s\colon W\to Z$ be the zero section. Let $p\colon W\to \{pt\}$. 
If $w\in \ec(Z)$ has the compact support, then we have
\[
p_*\circ\pi_*(w\cup s_*(1_Z))=p_*\circ s^*(w).
\]  
\end{enumerate}
\end{lem}
\begin{proof}
We can prove those claims by calculations of differential forms on the finite dimensional approximation $S^{2n+1}\times_{S^1}X$.
For (i), see Proposition VI.4.6 in \cite{audin}. 
The reader can refer to Proposition 6.15 and Proposition 6.24 in \cite{bott-tu} for non-equivariant version of (ii) and (iv).  
\end{proof}

\subsubsection{}\label{localization}

We set $\qf=\C(t)$ and $H^*_{S_1,\qf}(X)=H^*_{S_1}(X)\otimes_{\C[t]}\qf$.

Let $X$ be an $S^1$-manifold and $i$ denote the inclusion $\colon X^{S^1}\hookrightarrow X$. 
Let $\{Z_\lambda\}_{\lambda\in\Lambda}$ denote the set of connected components of the fixed point set, 
$i_\lambda$ denote the inclusion $Z_\lambda\hookrightarrow X$ and $\nu_\lambda$ denote the normal bundle on $Z_\lambda$ to $X$.
We can check that $\eu_{S^1}(\nu_\lambda)$ is invertible in $H^*_{S_1,\qf}(X)$.
\begin{thm}(\cite{atiyah-bott})
The map
\[
i_*\colon H^*_{S_1,\qf}(X^{S^1})\to H^*_{S_1,\qf}(X)
\]
is isomorphism.
In paticular, for $x\in H^*_{S_1,\qf}(X)$ we have
\[
x=\sum_{\lambda}{i_\lambda}_*\left(\frac{i_\lambda^*x}{\eu_{S^1}(\nu_\lambda)}\right).
\]
%where $e_{S^1}(\cdot)$ represents the $S^1$-equivariant Euler class.
\end{thm}

\subsection{Equivariant cohomology groups of quiver varieties}

%The arguments in \ref{spectral}--\ref{} work for arbitary quiver varieties. 

\subsubsection{}\label{vanish}
\begin{lem}
\[
H^m(\mathcal{M}_\zeta(\dimv))=0,\quad (\text{$m$ is odd,  or $m>n=\h{1}\mathrm{dim}\mathcal{M}_\zeta(\dimv)$}).
\]
\end{lem}
\begin{proof}
A quiver variety is endowed with a symplectic structure and an $S^1$-action which preserve the symplectic form. 
The moment map is a perfect Morse function (\cite{nakajima-lec-note} \S 5.1) and all the indices are even. 
So the odd degree cohomology groups vanish.

A quiver variety is homotopy equivalent to a certain Lagrangean subvariety (\cite{quiver1} Corollary 5.5). 
So the cohomology groups with degree larger than the half of the dimension vanish.
\end{proof}

\subsubsection{}\label{spectral}
\begin{prop}
There exists a (non-canonical) isomorphism as graded $\C[t]=H^*_{S^1}(\{pt\})$-module
\[
H^*_{S^1}(\mathcal{M}_\zeta(\dimv))\simeq H^*_{S^1}(\{pt\})\otimes H^*(\mathcal{M}_\zeta(\dimv)).
\]
\end{prop}
\begin{proof}
Since $BS^1$ is simply connected and the odd degree cohomologies of $BS^1$ and $\mathcal{M}_\zeta(\dimv)$ vanish, 
the spectral sequence associated with the fibration $ES^1\times_{S^1}\quiv{}{\dimv}\to BS^1$ degenerates at $E_2$-term. 
So the claim follows. 
\end{proof}

\begin{cor}
The forgetful map $H^*_{S^1}(\mathcal{M}_\zeta(\dimv))\to H^*(\mathcal{M}_\zeta(\dimv))$ is surjective.
\end{cor}
\begin{proof}
Note that the forgetful map can be described as $i^*$ 
where $i\colon \mathcal{M}_\zeta(\dimv)\hookrightarrow ES^1\times_{S^1}\mathcal{M}_\zeta(\dimv)$ is an inclusion of a fiber. 
Then the claim follows by Theorem 5.9 and Theorem 5.10 of \cite{guide-spectral}.
\end{proof}

\subsubsection{}

\begin{prop}
The following map is isomorphism :
\[
\begin{array}{ccc}
H^n_{S^1}(\quiv{}{\dimv})&\longrightarrow&H^{2n}_{S^1}(\quiv{}{\dimv})\\
x&\longmapsto&t^n\cdot x,
\end{array}
\]
where $n=\h{1}\mathrm{dim}(\quiv{}{\dimv})$.
\end{prop}
\begin{proof}
This follows from Lemma \ref{vanish} and Proposition \ref{spectral}.
\end{proof}

\subsubsection{}
Let us write $\ec(\quiv{}{\dimv})$ for $H^n_{S^1}(\quiv{}{\dimv})$, where $n=\h{1}\mathrm{dim}\,\quiv{}{\dimv}$.
\begin{dfn}
For $P\in(\quiv{}{\dimv})^{S^1}$ we define 
\[
%\xi_\lambda=t^{-n}\cdot i_{\lambda}^*(1_\lambda)
\xi_P=t^{-n}\cdot[P]\in \ec(\quiv{}{\dimv}).
\]
\end{dfn}

\subsubsection{}\label{basis}
For a finite dimensional $\C S^1$-module $M$, we define a number $\ueu(M)$ by the product of all the weight of $M$. 
Note that $\eu_{S^1}(M)=\ueu(M)\cdot t^{\mathrm{dim}M}$ when we regard $M$ as a $S^1$-equivariant vector bundle on a point.

For $P\in\quiv{}{\dimv}^{S^1}$ let $T_P$ denote the tangent space of $\quiv{}{\dimv}$ at $P$.

\begin{prop}
The set $\{\xi_P\}$ forms a basis of $\ec(\quiv{}{\dimv})$.
\end{prop}
\begin{proof}
Linear independency follows directly from Theorem \ref{localization}.
For $\alpha\in \ec(\quiv{}{\dimv})$ we have $i_P^*(\alpha)=c_P\cdot t^n$ for some $c_P\in\C$. So we have
\[
\alpha
=\sum_P \frac{c_P\cdot t^n\cdot{i_P}_*(1)}{\ueu(T_P)\cdot t^{2n}}
=\sum_P \frac{c_P}{\ueu(T_P)}\xi_P.
\]
Thus the claim follows. 
\end{proof}

\subsection{Representation of the affine Lie algebra}

\subsubsection{}\label{rep of asl}
Let $\mathbf{e}_i$ denote the $i$-th coordinate vector of $\Z^I$.
For $\dimv\in\Z^I$ we define the subvariety
\[
\B_i(\dimv)=\left\{(J_1,J_2)\in\quiv{0}{\dimv}\times\quiv{0}{\dimv+\mathbf{e}_i}\mid 
\text{$J_1$ is a subrepresentation of $J_2$}\right\}.
\]
of $\quiv{0}{\dimv}\times\quiv{0}{\dimv+\mathbf{e}_i}$. This is called {\bf Hecke correspondence}. 
This is smooth and $2|\dimv|+1$-dimensional (\cite{quiver2}).

Let $p_\varepsilon$ be a projection from $\quiv{0}{\dimv}\times\quiv{0}{\dimv+\mathbf{e}_i}$ to the $\varepsilon$-th factor.
We define operators $e_i$ and $f_i$ on $\oplus_{\dimv}H^{\mathrm{mid}}_{S^1}(\quiv{0}{\dimv})$ by
\begin{align*}
e_i(\alpha)&=(-1)^{v_{i-1}+v_{i}}{p_1}_*({p_2}^*(\alpha)\cup \B_i(\dimv))\quad\quad \left(\alpha\in{H}^{\mathrm{mid}}_{S^1}(\quiv{0}{\dimv+\mathbf{e}_i})\right),\\
f_i(\alpha)&=(-1)^{v_{i}+v_{i+1}}{p_2}_*({p_1}^*(\alpha)\cup \B_i(\dimv))\quad\quad \left(\alpha\in{H}^{\mathrm{mid}}_{S^1}(\quiv{0}{\dimv}\right).
\end{align*}
These operators give an representation of $\asl_l$.

\subsubsection{}
We define a bilinear form $\langle\ ,\ \rangle$ on $\lec(\quiv{0}{\dimv})$  by
\[
\langle \alpha,\beta\rangle=p_*(i_*)^{-1}(\alpha\cup \beta)\in H^*_{S^1,\qf}(\{pt\})
\]
where $i\colon \quiv{0}{\dimv}^{S^1}\hookrightarrow \quiv{0}{\dimv}$ and $p\colon \quiv{0}{\dimv}^{S^1}\to\{pt\}$.

For $\lambda\in\Pi$ we write simply $\xi_\lambda$ and $T_\lambda$ for $\xi_{J_\lambda}$ and $T_{J_\lambda}$.
%$T_{\lambda}$ denotes the tangent space of $\quiv{0}{\dimv}$ at $J_{\lambda}$ for $\lambda\in\Pi$.
For $\lambda$ and $\mu$ such that $(J_{\lambda},J_{\mu})\in\B_{i}(\dimv)$,
let $N_{\lambda,\mu}$ denote the fiber of the normal bundle on $\B_{i}(\dimv)$ 
to $\M(\dimv)\times\M(\dimv+\mathbf{e}_i)$ at $(J_{\lambda},J_{\mu})$.

\begin{prop}
\[
\langle \xi_{\lambda},\xi_{\mu}\rangle=\delta_{\lambda,\mu}\cdot\ueu(T_{\lambda}),\quad
\langle e_i\xi_{\lambda},\xi_{\mu}\rangle=(-1)^{v_{i-1}+v_{i}}\delta\bigl((J_\mu,J_\lambda)\in\B_i(\dimv)\bigr)\cdot\ueu(N_{\mu,\lambda}).
\]
\end{prop}
\begin{proof}
\begin{align*}
\langle \xi_{\lambda},\xi_{\mu}\rangle
&=t^{-|\lambda|-|\mu|} \left\langle {i_\lambda}_*(1),{i_\mu}_*(1_{\mu})\right\rangle\quad &\\
%&=t^{-|\lambda|-|\mu|} \left\langle 1_{\lambda},i_\lambda^*\left({i_\mu}_*(1_{\mu})\right)\right\rangle\quad &(\text{lemma \ref{lemma for ec} (i)})\\
&=t^{-|\lambda|-|\mu|} \left\langle 1,i_\lambda^*\circ{i_\mu}_*(1_{\mu})\right\rangle\quad &(\text{Lemma \ref{lemma for ec} (ii)})\\
&=t^{-|\lambda|-|\mu|} \left\langle 1,\delta_{\lambda,\mu}\cdot\eu(T_\lambda)\cdot1_{\mu}\right\rangle\quad &(\text{Lemma \ref{lemma for ec} (i)})\\
&=\ueu(T_\lambda).\quad&
\end{align*}
\begin{align*}
&(-1)^{v_{i-1}+v_{i}}\langle e_i\xi_{\lambda},\xi_{\mu}\rangle\\
&=t^{-|\lambda|-|\mu|}\left\langle {p_1}_*\bigl(p_2^*\circ{i_\lambda}_*(1)\cup [\B_i(\dimv)]\bigr),{i_\mu}_*(1)\right\rangle\quad&\\
&=t^{-|\lambda|-|\mu|}\left\langle {i_\mu}^*\circ{p_1}_*\bigl(p_2^*\circ{i_\lambda}_*(1)\cup [\B_i(\dimv)]\bigr),1\right\rangle\quad &(\text{Lemma \ref{lemma for ec} (ii)})\\
&=t^{-|\lambda|-|\mu|}\left\langle p_*\circ\bar{i_\mu}^*\bigl(p_2^*\circ{i_\lambda}_*(1)\cup [\B_i(\dimv)]\bigr),1\right\rangle\quad &(\text{Lemma \ref{lemma for ec} (iii)})\\
&=t^{-|\lambda|-|\mu|}\left\langle p_*\left({i_\lambda}_*(1)\cup \bar{i_\mu}^*\bigl([\B_i(\dimv)]\bigr)\right),1\right\rangle\quad&\\
&=t^{-|\lambda|-|\mu|}\left\langle {i_\lambda}^*\circ\bar{i_\mu}^*\bigl([\B_i(\dimv)]\bigr),1\right\rangle\quad &(\text{Lemma \ref{lemma for ec} (iv)})\\
&=t^{-|\lambda|-|\mu|}\langle\,\delta\bigl((J_\lambda,J_\mu)\in \B_i(\dimv)\bigr)\cdot\eu_{S^1}(N_{\mu,\lambda}),1\rangle\quad &(\text{Lemma \ref{lemma for ec} (i)})\\
&=\delta\bigl((J_\lambda,J_\mu)\in \B_i(\dimv)\bigr)\cdot\ueu(N_{\mu,\lambda}),
\end{align*}
where $p\colon X_2\to \{J_\mu\}$ and $\bar{i_\mu}\colon X_2\simeq \{J_\mu\}\times X_2\to X_1\times X_2$.
\end{proof}
\begin{cor}
\[
e_i\,\xi_{\lambda}=(-1)^{v_{i-1}+v_{i}}\sum_{\mu}\ueu(N_{\mu,\lambda}-T_{\lambda})\,\xi_{\mu}.
\]
where The summation runs over all $\mu$ such that $(J_{\mu},J_{\lambda})\in\B_{i}(\dimv)$, 
which is equivalent to that $\mu$ is obtained by removing a removable $i$-node from $\lambda$.
\end{cor}

\subsubsection{}\label{calofVV}
For a $\Z/l\Z$-module $M$, we set $M_i=\mathrm{Hom}_{\Z/l\Z}(\C_{(i)},M)$ where $\C_{(i)}$ is the 1-dimensional 
representation of $\Z/l\Z$ with weight $i$.

For two nodes $X=(a,b)$ and $X'=(a',b')$, we set $l(X,X')=a-b-a'+b'$.

\begin{prop}
\[
e_i\,\xi_\lambda=(-1)^{v_{i-1}+v_{i}}\sum_{\mu}\left(\prod_{A\in A_{\lambda,i}}-l(\lambda\backslash\mu,A)\prod_{R\in R_{\mu,i}}-l(\lambda\backslash\mu,R)^{-1}\right)\cdot\xi_\mu
\]
where the summation runs over all $\mu$ obtained by removing a removable $i$-node from $\lambda$.
\end{prop}
\begin{proof}
We write simply $V_\lambda$ for $V(J_\lambda)$. 
For a node $X=(a,b)$ we set $V_X=t^{a-b}$. 
Note that $V_\lambda=\sum_{X\in\lambda}V_X$.

Using the description in \cite{quiver2} (Corollary 3.12. and \S 5), we have
\begin{align*}
T_{\mu}&=\left((t+t^{-1}-2)V_{\mu}^*V_{\mu}+V_{\mu}+V_{\mu}^*\right)_0\\
N_{\mu,\lambda}&=\left((t+t^{-1}-2)V_{\mu}^*V_{\lambda}+V_{\lambda}+V_{\mu}^*-1\right)_0.
\end{align*}
Thus we have
\[
N_{\mu,\lambda}-T_{\lambda}=\left((t+t^{-1}-2)V_{\mu}^*+1\right)_i\cdot V_{\lambda\backslash\mu}-1.
\]
On the other hand, we can verify
\[
\left((t+t^{-1}-2)V_{\mu}+1\right)_i=\sum_{A\in A_{\mu,i}}\hspace{-4pt}V_A-\sum_{R\in R_{\mu,i}}\hspace{-4pt}V_R.
\]
Substituting this we get
\begin{align*}
N_{\mu,\lambda}-T_{\lambda}&=\sum_{A\in A_{\mu,i}}\hspace{-4pt}V_{\lambda\backslash\mu}V_A^*-\sum_{R\in R_{\mu,i}}\hspace{-4pt}V_{\lambda\backslash\mu}V_R^*-1\\
&=\sum_{A\in A_{\lambda,i}}\hspace{-4pt}V_{\lambda\backslash\mu}V_A^*-\sum_{R\in R_{\mu,i}}\hspace{-4pt}V_{\lambda\backslash\mu}V_R^*.
\end{align*}
Note that the weight of $V^{}_XV^*_{X'}$ equals to $l(X.X')$. 
So the claim follows.
\end{proof}

\subsubsection{}\label{thm2}

Let $L_\lambda$ be the product of all the negative weights of $T_{\lambda}$.
Normalize the basis by $b_\lambda=L_\lambda^{-1}\xi_\lambda$.
We identify $\oplus_{\dimv}H^{\mathrm{mid}}_{S^1}(\quiv{0}{\dimv})$ and $\B$ as vector spaces. 

\begin{thm}
The action on $\oplus_{\dimv}H^{\mathrm{mid}}_{S^1}(\quiv{0}{\dimv})$ coincides with representation defined in \ref{formulaforFrenkelKac}.
\end{thm}
\begin{proof}
We will check only for $e_i$'s. 

By \ref{tangentspace0} we have
$L_{\mu}=\Pi_{k,a}(-nl)$
where the product runs through all $k$ and $n>0$ such that $\maya_\mu(k)=-1$, $\maya_\mu(k+nl)=1$.
So we have
\begin{align*}
L_{\lambda}/L_{\mu}&=\prod_{A\in A_{\lambda,i}}-\big|l(\lambda\backslash\mu,A)\big|\prod_{R\in R_{\mu,i}}-\big|l(\lambda\backslash\mu,A)\big|^{-1}\\
&=(-1)^{\eta^+(\lambda\backslash\mu,i,\lambda)}\prod_{A\in A_{\lambda,i}}-l(\lambda\backslash\mu,A)\prod_{R\in R_{\mu,i}}-l(\lambda\backslash\mu,A).
\end{align*}
Note that 
\begin{align*}
\eta^+(\lambda\backslash\mu,i,\lambda)+\eta^-(\lambda\backslash\mu,i,\lambda)&\equiv c_{i-\h{1}}(\lambda)+c_{i+\h{1}}(\lambda)+1\\
&=v_{i-1}+v_{i}+v_{i}(\lambda)+v_{i+1}(\lambda).
\end{align*}
Thus the claim follows.
\end{proof}

\subsection{Representation of the Heisenberg algebra}\label{heisenberg}

\subsubsection{}\label{rep of heisen}
Recall that $M=\quiv{0}{\delta}$ is isomorphic to the minimal resolution of $\C^2/\hspace{-3pt}/(\Z/l\Z)$.
The exceptional fiber has $l-1$ connected component $\lag_i$ ($i=1,\ldots,l-1$)
and $\lag_i$ contains $P_{i-\h{1}}$ and $P_{i+\h{1}}$.

For $m\in\Z_{>0}$ we consider subvarieties
\[
\lag_i(m)=\{(I,J)\in M^{[n]}\times M^{[n+m]}\mid I\supset J,\ \mathrm{supp}(I/J)=\{x\} \ \text{for some $x\in\lag_i$}\}
\]
and define $\p_i(m)\colon\ec(M^{[n]})\to\ec(M^{[n+m]})$ by
\[
\p_i(m)(\alpha)=(-1)^m{p_2}_*(p_1^*(\alpha)\cup[\lag_i(m)])\quad (\alpha\in\ec(M^{[n]})).
\]
These operators  satisfy the relations of $\heisen_{l-1}$ (\cite{nakajima-heisen}).

\subsubsection{}\label{thmofQW}
%Let $i_{\blam}$ denote the inclusion $\{J_{\blam}\}\hookrightarrow\quiv{\infty}{n\delta}$ 
%and $1_{\blam}$ denote the generator of $H_{S^1}^0(\{J_{\blam}\})$.
%As befor we can define 
%$\xi_\blam=t^{-n}[J_\blam]\in H^{\mathrm{mid}}_{S^1}(\quiv{\infty}{n\delta})$ 
%and they form a basis. 

For $\blam\in\Pi^{\tilde{I}}$ let $L_\blam$ denote the product of all the negative weights of $T_\blam$.
%let $N_{\blam}$ denote the products of all the negative weights of the tangent space of $\quiv{\infty}{n\delta}$ at $\blam$. 
We set $b_{\blam}=L_{\blam}^{-1}\xi_{J_\blam}$. 
Let us identify $\oplus_n H^{\mathrm{mid}}_{S^1}(M^{[n]})$ and $\B^{\otimes \tilde{I}}$ as vector spaces.

\begin{prop}(\cite{qin-wang} Lemma 3.3)
The action of $\heisen_{l-1}$ on $\oplus_n H^{\mathrm{mid}}_{S^1}(M^{[n]})$ coincides with the action defined in \ref{higherbosonrep}.
\end{prop}

\subsection{Main theorem for equivariant cohomologies}\label{main for ec}
The $S^1$-equivariant diffeomorphism induces the isomorphism 
\[
\oplus_{\dimv}H^{\mathrm{mid}}_{S^1}(\quiv{0}{\dimv})\simeq \oplus_{\dimv}H^{\mathrm{mid}}_{S^1}(\quiv{\infty}{\dimv}).
\]
The quiver variety $\quiv{\infty}{\dimv}$ is isomorphic to the Hilbert scheme and we have the isomorphism
\[
\oplus_\dimv \phi_{\dimv}^*\colon \oplus_{\dimv}H^{\mathrm{mid}}_{S^1}(\quiv{\infty}{\dimv})\simeq \C Q\otimes \left(\oplus_n H^{\mathrm{mid}}_{S^1}(M^{[n]})\right).
\]
The $\asl_l$-action on $\oplus_{\dimv}H^{\mathrm{mid}}_{S^1}(\quiv{0}{\dimv})$ is given in \ref{rep of asl}
Apply the Frenkel-Kac construction for the $\heisen_{l-1}$-action on $\oplus_n H^{\mathrm{mid}}_{S^1}(M^{[n]})$ given in \ref{rep of heisen}, then we have the $\asl_l$-action on $\C Q\otimes \left(\oplus_n H^{\mathrm{mid}}_{S^1}(M^{[n]})\right)$.

\begin{thm}
The composition of the two isomorphism 
\[
\oplus_{\dimv}H^{\mathrm{mid}}_{S^1}(\quiv{0}{\dimv})\simeq \oplus_{\dimv}H^{\mathrm{mid}}_{S^1}(\quiv{\infty}{\dimv})\simeq \C Q\otimes \left(\oplus_n H^{\mathrm{mid}}_{S^1}(M^{[n]})\right).
\]
intertwines the $\asl_l$-actions.
\end{thm}
\begin{proof}
This follows from Theorem \ref{formulaforFrenkelKac}, Theorem \ref{thm1}, Theorem \ref{thm2} and Proposition \ref{thmofQW}.
\end{proof}

\section{Representations on ordinary cohomologies}\label{ordinary}

\subsection{Representations on ordinary cohomologies}
The constructions of representations of the affine Lie algebra and the Heisenberg algebra on the equivariant cohomology groups 
in \ref{rep of asl} and \ref{rep of heisen} can be applied to the ordinary cohomology groups too.

\begin{thm}(\cite{quiver2})
The affine Lie algebra $\asl_l$ acts on $\oplus_{\dimv}H^{\mathrm{mid}}(\quiv{0}{\dimv})$, 
and this is a level-$1$ integrable highest weight representation. 
%In particular this is irreducible.  
\end{thm}

\begin{thm}(\cite{nakajima-heisen})
The Heisenberg algebra $\heisen_{l-1}$ acts on $\oplus_{n}H^{\mathrm{mid}}(\quiv{0}{n\delta})$, 
and this is the Fock space representation. 
%In particular this is irreducible.  
\end{thm}

\subsection{Main theorem for ordinary cohomologies}\label{main for ord}
As the case of equivariant cohomology groups in \S \ref{main for ec}, we have the canonical isomorphisms of ordinary cohomology groups 
%\[
%\oplus_{\dimv}H^{\mathrm{mid}}(\quiv{0}{\dimv})\simeq \oplus_{\dimv}H^{\mathrm{mid}}(\quiv{\infty}{\dimv}),
%\]
%\[
%\oplus_\dimv \phi_{c,n}^*\colon \oplus_{\dimv}H^{\mathrm{mid}}(\quiv{\infty}{\dimv})\simeq \C Q\otimes \left(\oplus_n H^{\mathrm{mid}}(M^{[n]})\right),
%\]
and the $\asl_l$-actions on $\oplus_{\dimv}H^{\mathrm{mid}}(\quiv{0}{\dimv})$ and $\C Q\otimes \left(\oplus_n H^{\mathrm{mid}}(M^{[n]})\right)$.

\begin{thm}
The composition of the two isomorphisms 
\[
\oplus_{\dimv}H^{\mathrm{mid}}(\quiv{0}{\dimv})\simeq \oplus_{\dimv}H^{\mathrm{mid}}(\quiv{\infty}{\dimv})\simeq \C Q\otimes \left(\oplus_n H^{\mathrm{mid}}(M^{[n]})\right).
\]
intertwines the $\asl_l$-actions.
\end{thm}
\begin{proof}
Note that the forgetful map is compatible with the actions and the isomorphisms. 
Then the claim follows from Corollary \ref{spectral} and Theorem \ref{main for ec}.
\end{proof}

\providecommand{\bysame}{\leavevmode\hbox to3em{\hrulefill}\thinspace}
\providecommand{\MR}{\relax\ifhmode\unskip\space\fi MR }
\providecommand{\MRhref}[2]{%
  \href{http://www.ams.org/mathscinet-getitem?mr=#1}{#2}
}
\providecommand{\href}[2]{#2}

{\tt
\noindent Kentaro Nagao

\noindent Department of Mathematics, Kyoto University, Kyoto 606-8502, Japan

\noindent kentaron@math.kyoto-u.ac.jp
}
\end{document}